\documentclass[a4paper,oneside,11pt]{amsart}

\usepackage[dvipdfmx]{graphicx} 
\usepackage{pdfpages}
\usepackage{float}
\usepackage{amssymb}

\usepackage{bm}

\newtheorem{lem}{Lemma}[section]
\newtheorem{prop}[lem]{Proposition}
\newtheorem{thm}[lem]{Theorem}
\newtheorem{cor}[lem]{Corollary}

\newtheorem{df}[lem]{Definition}
\newtheorem{rem}[lem]{Remark}
\newtheorem{exmp}[lem]{Example}

\newcommand{\vertiii}[1]{{\left\vert\kern-0.25ex\left\vert\kern-0.25ex\left\vert #1
    \right\vert\kern-0.25ex\right\vert\kern-0.25ex\right\vert}}

\numberwithin{equation}{section}

\allowdisplaybreaks

\title[Dynamics of the multicolor box-ball system]{Dynamics of the multicolor box-ball system \\ with random initial conditions\\ via Pitman's transformation}

\author{Kazuki Kondo}

\begin{document}

\maketitle
\setcounter{tocdepth}{3}

\tableofcontents

\section{Introduction}

The Box-Ball System (BBS) is a one-dimensional cellular automaton in $\{0,1\}^{\mathbb{Z}}$ that was introduced by Takahashi and Satsuma in 1990 \cite{TS}, and has been extensively studied from the viewpoint of integrable systems. In particular, it is connected with the KdV equation \cite{KdV} 
\[\frac{\partial u}{\partial t}+6u\frac{\partial u}{\partial x}+\frac{\partial^3 u}{\partial x^3}=0,\ \ u=u(x,t),\ x,t\in\mathbb{R},\]
which is a non-linear partial differential equation giving a mathematical model for waves on shallow water surfaces. The BBS equation of motion is obtained from the KdV equation by applying an appropriate discretization and transform \cite{TTMS}. The KdV equation has soliton solutions whose shape and speed are conserved after collision with other solitons, and such a phenomenon is also observed in the BBS. 

Now we present the original definition of the BBS from \cite{TS}. We denote a particle configuration by $(\eta_n)_{n\in\mathbb{Z}}\in \{0,1\}^{\mathbb{Z}}$ for the two-sided case or $(\eta_n)_{n\in\mathbb{N}}\in \{0,1\}^{\mathbb{N}}$ for the one-sided case. Specifically, we write $\eta_n = 1$ if there is a particle at site $n$, and $\eta_n = 0$ otherwise. On the condition that there is a finite number of particles, that is, $\sum_{n\in\mathbb{Z}}\eta_n<\infty$, the evolution of the BBS is described by an operator $T:\{0,1\}^{\mathbb{Z}}\rightarrow\{0,1\}^{\mathbb{Z}}$ that is characterized by the following BBS equation of motion, 
\[(T\eta)_{n}=\min\left\{1-\eta_{n},\sum_{m=-\infty}^{n-1}\left(\eta_m - (T\eta)_m\right)\right\},\]
where we suppose $(T\eta)_n=0$ for $n\leq \inf\{l:\eta_l=1\}$, so the sums in the above definition are well-defined. In other words, the balls move sequentially from left to right, that is, from negative to positive, with each being transported to the leftmost unoccupied site to its right as follows.

\vspace{10pt}

\makebox[2.6cm][r]{$\eta=\ $}$(\cdots\ 0\ 1\ 1\ 1\ 0\ 0\ 0\ 0\ 0\ 0\ 1\ 0\ 0\ 0\ 0\ 0\ 0\ 0\ 0\ 0\ 0\ \cdots)$

\makebox[2.6cm][r]{$T\eta=\ $}$(\cdots\ 0\ 0\ 0\ 0\ 1\ 1\ 1\ 0\ 0\ 0\ 0\ 1\ 0\ 0\ 0\ 0\ 0\ 0\ 0\ 0\ 0\ \cdots)$

\makebox[2.6cm][r]{$T^2\eta=\ $}$(\cdots\ 0\ 0\ 0\ 0\ 0\ 0\ 0\ 1\ 1\ 1\ 0\ 0\ 1\ 0\ 0\ 0\ 0\ 0\ 0\ 0\ 0\ \cdots)$

\makebox[2.6cm][r]{$T^3\eta=\ $}$(\cdots\ 0\ 0\ 0\ 0\ 0\ 0\ 0\ 0\ 0\ 0\ 1\ 1\ 0\ 1\ 1\ 0\ 0\ 0\ 0\ 0\_0\ \cdots)$

\makebox[2.6cm][r]{$T^4\eta=\ $}$(\cdots\ 0\ 0\ 0\ 0\ 0\ 0\ 0\ 0\ 0\ 0\ 0\ 0\ 1\ 0\ 0\ 1\ 1\ 1\ 0\ 0\ 0\ \cdots)$
\vspace{10pt}

\hspace{-13pt}This example exhibits a string of $3$ consecutive balls, called a soliton, moving distance $3$ in each time step when there is no interaction, and recovering its shape and speed after a collision with another soliton (of length $1$).

In this paper we consider a generalization of the BBS that incorporates multiple colors of balls, that is, we assume that there are $\kappa$-color balls (particles) for some $\kappa\in\mathbb{N}$. 
This model is called multicolor BBS and was introduced in \cite{TaK}, as a generalization of the original $\kappa=1$ BBS first introduced in \cite{T}. In this model, particle configurations are given by $(\eta_n)_{n\in\mathbb{Z}}\in \{0,1,\cdots,\kappa\}^{\mathbb{Z}}$, where we suppose that the numbers $1,\cdots,\kappa$ represent the colors of the balls and $0$ represents the empty site. For each $i=1,\cdots,\kappa$, we define the operator $T_i$ under which the balls of color $i$ move from left to right, with each being transported to the leftmost unoccupied site to its right, with balls of other colors remaining static. The dynamics of the multicolor BBS are then defined by the operator $T=T_\kappa\circ\cdots\circ T_1$.

For example, the evolution of the BBS with 3-color balls is as follows 

\vspace{10pt}

\makebox[2.6cm][r]{$\eta=\ $}$(\cdots\ 0\ 1\ 2\ 0\ 3\ 1\ 3\ 2\ 0\ 3\ 0\ 1\ 1\ 2\ 3\ 0\ 0\ 0\ 0\ 0\ 0\ 0\ 0\ 0\ \cdots)$

\makebox[2.6cm][r]{$T_1\eta=\ $}$(\cdots\ 0\ 0\ 2\ 1\ 3\ 0\ 3\ 2\ 1\ 3\ 0\ 0\ 0\ 2\ 3\ 1\ 1\ 0\ 0\ 0\ 0\ 0\ 0\ 0\ \cdots)$

\makebox[2.6cm][r]{$T_2\circ T_1\eta=\ $}$(\cdots\ 0\ 0\ 0\ 1\ 3\ 2\ 3\ 0\ 1\ 3\ 2\ 0\ 0\ 0\ 3\ 1\ 1\ 2\ 0\ 0\ 0\ 0\ 0\ 0\ \cdots)$

\makebox[2.6cm][r]{$T\eta=T_3\circ T_2\circ T_1\eta=\ $}$(\cdots\ 0\ 0\ 0\ 1\ 0\ 2\ 0\ 3\ 1\ 0\ 2\ 3\ 3\ 0\ 0\ 1\ 1\ 2\ 3\ 0\ 0\ 0\ 0\ 0\ \cdots)$

\makebox[2.6cm][r]{$T^2\eta=\ $}$(\cdots\ 0\ 0\ 0\ 0\ 1\ 0\ 2\ 0\ 3\ 1\ 0\ 0\ 0\ 2\ 3\ 3\ 0\ 0\ 0\ 1\ 1\ 2\ 3\ 0\ \cdots)$

\vspace{10pt}

\hspace{-13pt}where $T=T_3\circ T_2\circ T_1$. In the multicolor case, a string of consecutive balls of non-decreasing colors is called a soliton and shows the same behavior as in the 1-color case.

The multicolor BBS with finite number of balls has been well studied mostly in the context of integrable systems (see, e.g., the review article \cite{IKT} or the textbook on the BBS \cite{Tokihiro}). Recently, \cite{KLO} and \cite{AH}, \cite{LLPS} considered the multicolor BBS with one-sided random initial configuration and derived scaling limits of probability measures on the space of κ-tuple of Young diagrams induced by the random configuration. Later, we introduce the two-sided version of the multicolor BBS, which is one of the main contributions of this paper.

The dynamics of the one-color BBS has been extended to two-sided infinite configurations and studied when the initial condition is random \cite{DS,F}. In the paper \cite{DS} for the one-color BBS, the particle configuration is encoded by a certain path $S=(S_n)_{n\in\mathbb{Z}}$ in $\mathbb{Z}$ and the action $T$ of the BBS is defined via an operation on the path space. Moreover, a formal inverse $T^{-1}$ of $T$ is defined, and the class of configurations $S$ below such that $TS$ and $T^{-1}S$ are well-defined and reversible for all times, i.e. 
\[\{S=(S_n)_{n\in\mathbb{Z}} \in \mathbb{R}^{\mathbb{Z}}\::\:T^{k}S \mbox{\ is well defined and }TT^{-1}(T^{k}S)=T^{-1}T(T^{k}S)=T^{k}S,\ \forall k\in\mathbb{Z}\} ,\]
is precisely characterized. Within this framework, random initial conditions such that almost all paths are in the class is studied from the viewpoint of invariance under $T$, the current of particles crossing the origin, and the speed of a single tagged particle. 

Such an extended analysis was made possible thanks to connection that was identified between the BBS dynamics and Pitman's transformation. Indeed, in \cite{DS}, the action $T$ on the path space is shown to correspond to the operation of reflection in the past maximum of the path, which is precisely the operation known Pitman transform. Pitman transform is introduced by \cite{Pitman} and appears in the well-known Pitman's theorem, which states that if $(B_t)_{t\ge0}$ is a one-dimensional Brownian motion, then the stochastic process 
$(2\sup_{0\le s\le t}B_s-B_t)_{t\ge0}$ is a three dimensional Bessel process, i.e. is distributed as the Euclidean norm of a three dimensional Brownian motion. This transform has been generalized to the multidimensional case by Biane \cite{BBC}, and in this paper, we show that the actions of the multicolor BBS can be described by the multidimensional Pitman transform.

We start by introducing the one-sided and two-sided Pitman transform for the multicolor BBS theory (Section \ref{one-sided pitman}, \ref{two-sided pitman and its inverse}). Next, as in the case of the one-color BBS, we show that particle configurations of multicolor BBS can be encoded by a certain path in $\mathbb{R}^{\mathbb{\kappa}}$ (Section \ref{Vectors for path encodings}, \ref{one-sided configuration}) and the action $T_i$ corresponds to the composition of the extended Pitman transform and a certain operator (Section \ref{one-sided BBS}, \ref{two-sided BBS}). Moreover, we characterize the set of configurations for which the actions $T_1,T_2,\cdots T_\kappa$ are well-defined and reversible for all times (Section \ref{invariant set}). Then, we give an example of a random initial condition that is invariant in distribution under the dynamics of the multicolor BBS (Section \ref{iid}). Finally, we consider a generalization of the multicolor BBS, that is defined for continuous paths on $\mathbb{R}$ (Section \ref{onR}), and show that $\kappa$-dimensional Brownian motion with a certain drift is invariant under the action of the generalized multicolor BBS (Section \ref{BMD}).

Regarding notational conventions, we distinguish $\mathbb{N}=\{1,2,\dots,\}$ and $\mathbb{Z}_+=\{0,1,\dots\}$. 

\section{Pitman transform}\label{pitman}

In this section, we prepare Pitman transform and the extended versions of it which will be used for the path encoding of the particle configuration in the subsequent sections. We start by defining one-sided Pitman transform and studying its property (Section \ref{one-sided pitman}). Then, in Section \ref{two-sided pitman and its inverse}, we define two-sided Pitman transform and examine its inverse on an appropriate set.

\subsection{One-sided Pitman transform}\label{one-sided pitman}
We first see the definition of the multidimensional version of Pitman transform introduced by Biane \cite{BBC}.
{\df \label{original pitman}
Suppose that ${\mathbb{R}}^k$ is k-dimensional Euclidean space with dual space $V$ and let $\alpha\in {\mathbb{R}}^k, \alpha^*\in V$ be such that $\alpha^*(\alpha)=2$. The Pitman transform  $P_\alpha$ is defined on the set of continuous paths $\pi:[0,T]\to {\mathbb{R}}^k$, satisfying $\pi(0)=0$, by the formula, 
\[P_{\alpha,\alpha^*}\pi(t)=\pi(t)-\inf_{0\leq s\leq t}\alpha^*(\pi(s))\alpha,\ \ \ 0\leq t\leq T. \]}

For the multicolor BBS theory, we take the domain of $\pi$ as $\mathbb{Z}_+$ and $\alpha^{*}$ as the inner product with $\frac{\alpha}{|\alpha|^2}$ in the above definition, and define the one-sided Pitman transform.

{\df \label{oneside pitman}
Let $\alpha\in{\mathbb{R}}^k$ be such that $\alpha\ne0$. The one-sided Pitman transform with respect to $\alpha$ is defined on the set of discrete paths $\pi:\mathbb{Z_+}\to {\mathbb{R}}^k$, satisfying $\pi(0)=0$, by the formula, 
\[P_\alpha\pi(n)=\pi(n)-2\inf_{0\leq m\leq n}\frac{\alpha\cdot\pi(m)}{|\alpha|^2}\alpha,\ \ \ n\geq 0 \]
where $\alpha\cdot\pi(m)$ is the inner product of $\alpha$ and $\pi(m)$, and $|\alpha|^2=\alpha\cdot\alpha$.}

{\exmp \label{1dimpitman} For any $\alpha\in\mathbb{R},\ \alpha\ne0$, the one-sided Pitman transform is given by\[P_\alpha\pi(n)=\pi(n)-2\inf_{0\leq m\leq n}\pi(m),\ \ \ n\geq 0 \ \]for $\pi:\mathbb{Z_+}\to {\mathbb{R}}$, satisfying $\pi(0)=0$. Therefore the one-sided Pitman transform $P_\alpha$ on 1-dimensional Euclidean space does not depend on $\alpha$. We write it as  $P_1$. (See Figure 1.)}

{\df \label{1dim oneside pitman}\[P_1:=P_\alpha\ \ \ for\ \alpha \in\mathbb{R},\ \alpha\ne0.\]
That is, 
\begin{equation}\label{half P_1}
P_1\pi(n)=\pi(n)-2\inf_{0\leq m\leq n}\pi(m),\ \ \ n\geq 0
\end{equation}
for $\pi:\mathbb{Z_+}\to {\mathbb{R}}$, satisfying $\pi(0)=0$.}

\begin{figure}[H]
\centering
\scalebox{0.40}{\includegraphics{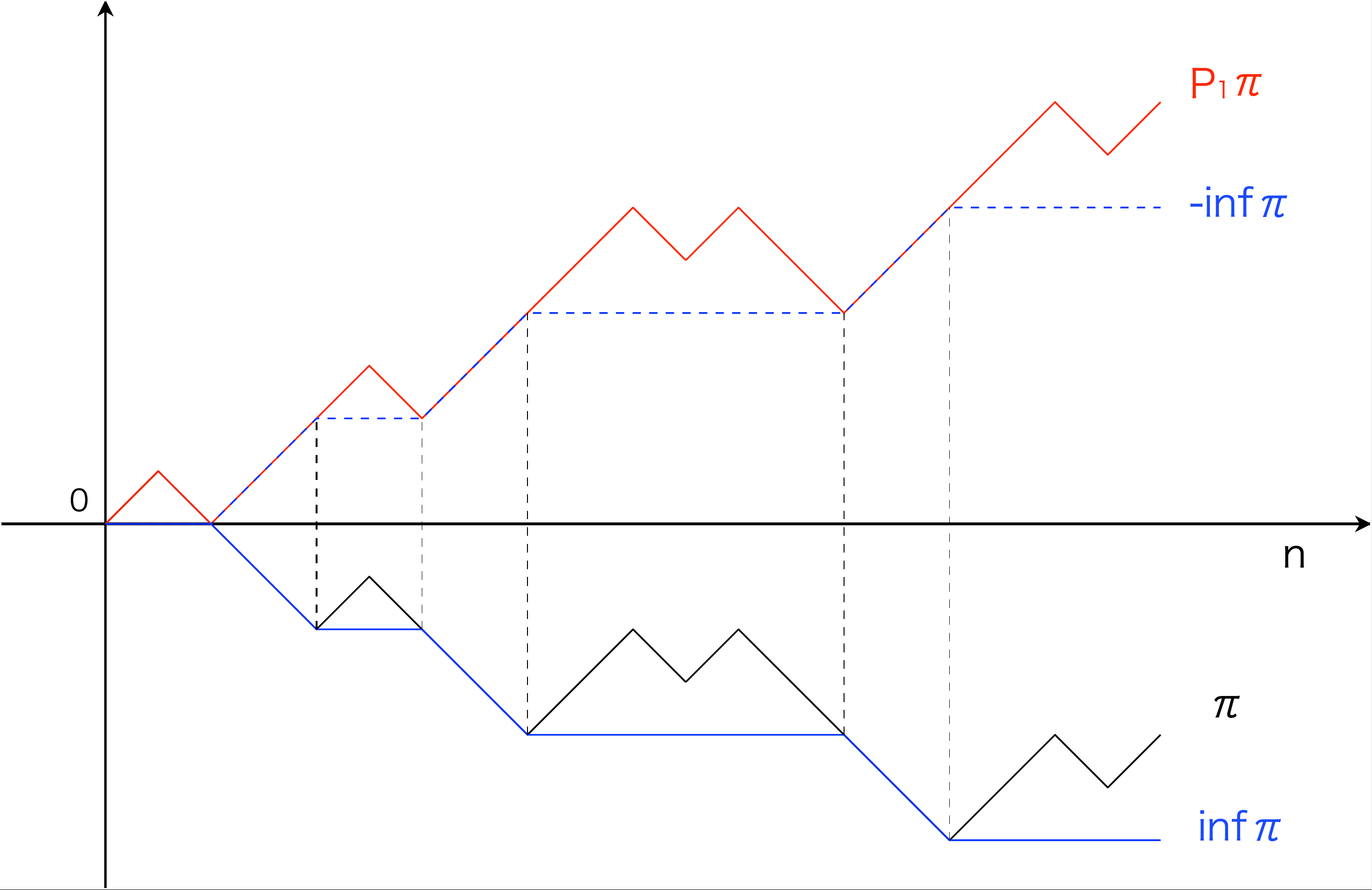}}
\vspace{10pt}
\caption{$P_1\pi(n)=\pi(n)-2\inf_{0\leq m\leq n}\pi(m).$}
\end{figure}

Next, we show the useful property of the one-sided Pitman transform for considering the actions of the BBS.

{\prop\label{half orthogonal} Let $k\geq2$ and $\pi_\alpha(n):=\frac{\alpha\cdot\pi(n)}{|\alpha|^2}$ for $\alpha\in\mathbb{R}^k$. $\pi:\mathbb{Z_+}\to {\mathbb{R}}$ is decomposed into the sum of the vector projection of $\pi(n)$ along $\alpha$ and the vector orthogonal to $\alpha$\::
\[\pi(n)=\pi_\alpha(n)\alpha+\left\{\pi(n)-\pi_\alpha(n)\alpha\right\}\]
for any $n\geq0$. Then, it holds that
\[P_\alpha\pi(n)=\left\{P_1\pi_\alpha(n)\right\}\alpha+\left\{\pi(n)-\pi_\alpha(n)\alpha\right\}.\](See figure 2.)}
\begin{proof}
\begin{align*}
P_\alpha\pi(n)&=\pi(n)-2\inf_{0\leq m\leq n}\frac{\alpha\cdot\pi(m)}{|\alpha|^2}\alpha\\
&=\pi_\alpha(n)\alpha+\left\{\pi(n)-\pi_\alpha(n)\alpha\right\}-2\inf_{0\leq m\leq n}\pi_\alpha(m)\alpha\\
&=\left\{\pi_\alpha(n)-2\inf_{0\leq m\leq n}\pi_\alpha(m)\right\}\alpha+\left\{\pi(n)-\pi_\alpha(n)\alpha\right\}\\
&=\left\{P_1\pi_\alpha(n)\right\}\alpha+\left\{\pi(n)-\pi_\alpha(n)\alpha\right\}.
\end{align*}
\end{proof}

\begin{figure}[H]
\centering
\scalebox{0.40}{\includegraphics{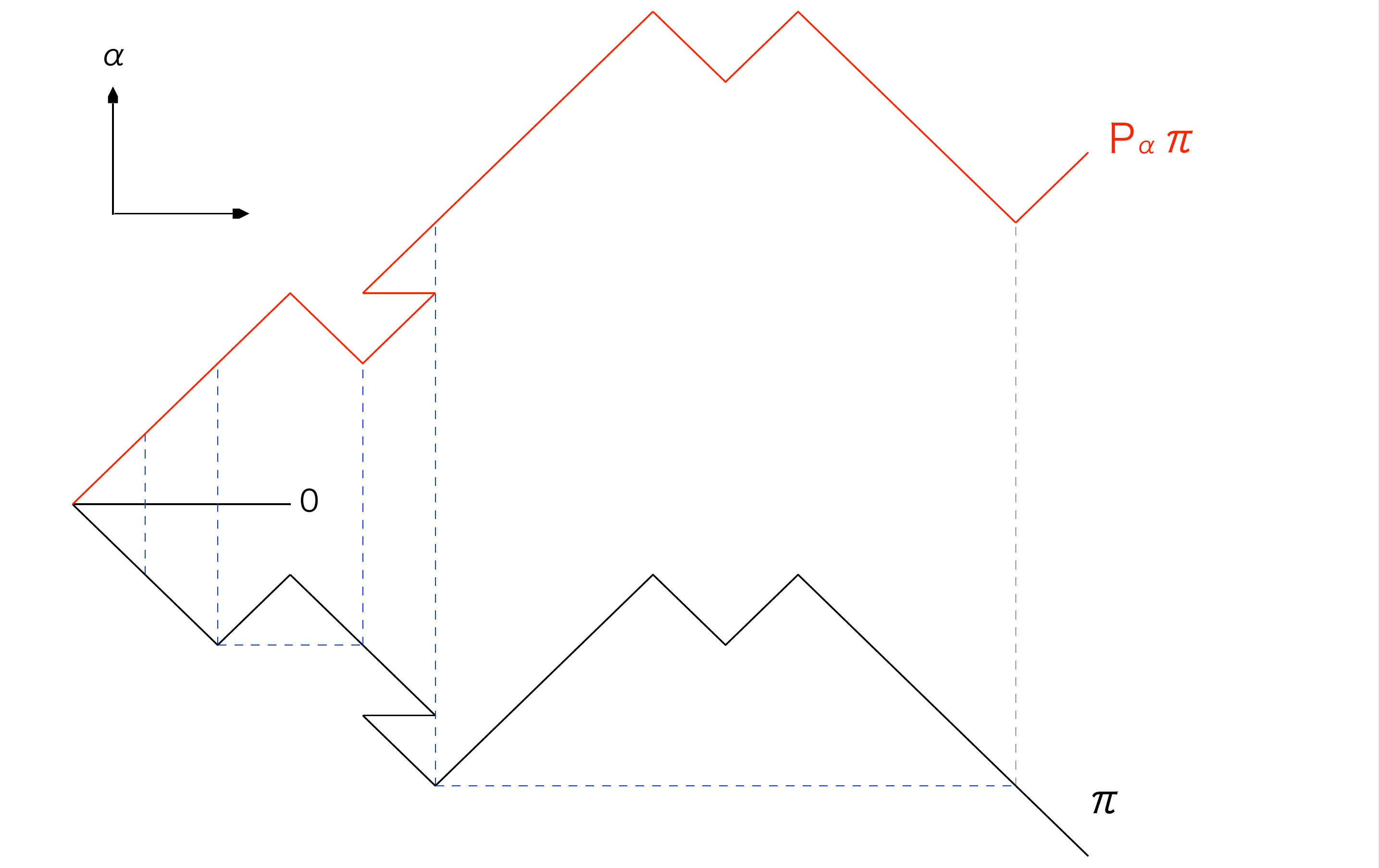}}
\vspace{0pt}
\caption{$P_\alpha\pi(n)=\left\{P_1\pi_\alpha(n)\right\}\alpha+\left\{\pi(n)-\pi_\alpha(n)\alpha\right\}.$}
\end{figure}

\subsection{Two-sided Pitman transform and its inverse}\label{two-sided pitman and its inverse}

This section provides the two-sided Pitman transform and its inverse on an appropriate set.

{\df \label{twoside pitman}
Let $\alpha\in{\mathbb{R}}^k,\ \alpha\ne0$. The two-sided Pitman transform with respect to $\alpha$ is defined on the set of discrete paths \[\{\pi:\mathbb{Z}\to {\mathbb{R}}^k,\ \pi(0)=0,\ \inf_{m\leq0}\alpha\cdot\pi(m)>-\infty\}\]by the formula, 
\[P_\alpha\pi(n)=\pi(n)-2\inf_{m\leq n}\frac{\alpha\cdot\pi(m)}{|\alpha|^2}\alpha+2\inf_{m\leq 0}\frac{\alpha\cdot\pi(m)}{|\alpha|^2}\alpha,\ \ \ n\in\mathbb{Z} .\]
Similarly to Example \ref{1dimpitman}, it holds that 

\[P_\alpha\pi(n)=\pi(n)-2\inf_{m\leq n}\pi(m)+2\inf_{m\leq 0}\pi(m),\ \ \ n\in\mathbb{Z}\]
for any $\alpha\in\mathbb{R},\ \alpha\ne0$, and it does not depend on $\alpha$. Then we define

\[P_1:=P_\alpha\ \ \ for\ \alpha \in\mathbb{R},\ \alpha\ne0.\]
That is, 
\begin{equation}\label{P_1}
P_1\pi(n)=\pi(n)-2\inf_{m\leq n}\pi(m)+2\inf_{m\leq 0}\pi(m),\ \ \ n\in\mathbb{Z} \ 
\end{equation}
for $\pi:\mathbb{Z}\to {\mathbb{R}}$, satisfying $\pi(0)=0,\ \inf_{m\leq 0}\pi(m)>-\infty$.

}

Next, we introduce a new transform which will be inverse of the two-sided Pitman transform on an appropriate set.

{\df \label{inverse pitman}
Let $\alpha\in{\mathbb{R}}^k,\ \alpha\ne0$. Define the transform ${P_\alpha}^{-1}$ on the set of discrete paths \[\{\pi:\mathbb{Z}\to {\mathbb{R}}^k,\ \pi(0)=0,\ \inf_{m\geq0}\alpha\cdot\pi(m)>-\infty\}\]by the formula, 
\[{P_\alpha}^{-1}\pi(n)=\pi(n)-2\inf_{m\geq n}\frac{\alpha\cdot\pi(m)}{|\alpha|^2}\alpha+2\inf_{m\geq 0}\frac{\alpha\cdot\pi(m)}{|\alpha|^2}\alpha,\ \ \ n\in\mathbb{Z} .\]
In this case, it also holds that 
\[P^{-1}_\alpha\pi(n)=\pi(n)-2\inf_{m\geq n}\pi(m)+2\inf_{m\geq 0}\pi(m),\ \ \ n\in\mathbb{Z}\]
for any $\alpha\in\mathbb{R},\ \alpha\ne0$, and it does not depend on $\alpha$. Then we define
\[{P_1}^{-1}:=P^{-1}_\alpha\ \ \ for\ \alpha \in\mathbb{R},\ \alpha\ne0.\]
That is, 
\begin{equation}\label{P^-1_1}
{P_1}^{-1}\pi(n)=\pi(n)-2\inf_{m\geq n}\pi(m)+2\inf_{m\geq 0}\pi(m),\ \ \ n\in\mathbb{Z} \ 
\end{equation}
for $\pi:\mathbb{Z}\to {\mathbb{R}}$, satisfying $\pi(0)=0,\ \inf_{m\geq0}\pi(m)>-\infty$.}

{\rem \label{orthogonal}With the same notation as Proposition \ref{half orthogonal}, it holds that,
\[P_\alpha\pi(n)=\left\{P_1\pi_\alpha(n)\right\}\alpha+\left\{\pi(n)-\pi_\alpha(n)\alpha\right\}\]
\[P^{-1}_\alpha\pi(n)=\left\{P^{-1}_1\pi_\alpha(n)\right\}\alpha+\left\{\pi(n)-\pi_\alpha(n)\alpha\right\}.\]
Therefore, $P_1 {P_1}^{-1}=\mathrm{id.}$ on some set $E_\alpha$ implies $P_\alpha {P_\alpha}^{-1}=\mathrm{id.}$ on $\{\pi :\mathbb{Z}\to {\mathbb{R}}^k\::\:\pi_\alpha\in\ E_\alpha\}$, and ${P_1}^{-1}P_1=\mathrm{id.}$ on some set $F_\alpha$ implies $ {P_\alpha}^{-1}P_\alpha=\mathrm{id.}$ on $\{\pi :\mathbb{Z}\to {\mathbb{R}}^k\::\:\pi_\alpha\in\ F_\alpha\}$.
}

{\df We define the domain of $P_1$ and ${P_1}^{-1}$, and their subsets,
\begin{equation}\label{P}
\mathcal{R}^{P_1}:=\{\pi:\mathbb{Z}\to {\mathbb{R}},\ \ \pi(0)=0,\ \inf_{m\leq0}\pi(m)>-\infty\},
\end{equation}

\begin{equation}\label{P'}
\mathcal{R}^{{P_1}^{-1}}:=\{\pi:\mathbb{Z}\to {\mathbb{R}},\ \ \pi(0)=0,\ \inf_{m\geq0}\pi(m)>-\infty\},
\end{equation}

\begin{equation}\label{P'P}
\mathcal{R}^{{P_1}^{-1}P_1}:=\{\pi\in\mathcal{R}^{P_1}\::\:|\pi(n+1)-\pi(n)|\in\{0,1\},\ \forall n,\ \inf_{m\leq n}\pi(m)=\pi(n)\ \ i.o.\ as\ \ n\rightarrow\infty\},
\end{equation}

\begin{equation}\label{PP'}
\mathcal{R}^{P_1{P_1}^{-1}}:=\{\pi\in\mathcal{R}^{{P_1}^{-1}}\::\:|\pi(n+1)-\pi(n)|\in\{0,1\},\ \forall n,\ \inf_{m\geq n}\pi(m)=\pi(n)\ \ i.o.\ as\ \ n\rightarrow-\infty\}.
\end{equation}}

We prepare following proposition to guarantee that ${P_1}^{-1}P_1$ and $P_1{P_1}^{-1}$ are well-defined on $\mathcal{R}^{P_1}$ and $\mathcal{R}^{{P_1}^{-1}}$ respectively.

{\prop It holds that
\[P_1\left(\mathcal{R}^{P_1}\right)\subseteq\mathcal{R}^{{P_1}^{-1}},\]
\[{P_1}^{-1}\left(\mathcal{R}^{{P_1}^{-1}}\right)\subseteq\mathcal{R}^{P_1}.\]}

\begin{proof}
Suppose that $n\geq0$ and $\pi\in\mathcal{R}^{P_1}$. Since $\inf_{m\leq n}\pi(m)\leq \inf_{m\leq 0}\pi(m)$, we have

\begin{align*}
P_1\pi(n)&=\pi(n)-2\inf_{m\leq n}\pi(m)+2\inf_{m\leq 0}\pi(m)\\
&\geq \pi(n).
\end{align*}
On the other hand, since $\inf_{m\leq n}\pi(m)\leq \pi(n)$, we have
\begin{align*}
P_1\pi(n)&=\pi(n)-2\inf_{m\leq n}\pi(m)+2\inf_{m\leq 0}\pi(m)\\
&\geq -\pi(n)+2\inf_{m\leq 0}\pi(m).
\end{align*}
The above two inequalities show 
\[P_1\pi(n)-\inf_{m\leq 0}\pi(m)\geq \pm\left\{-\pi(n)+\inf_{m\leq 0}\pi(m)\right\},\]
then
\[P_1\pi(n)\geq\inf_{m\leq 0}\pi(m).\]
It shows the first claim and we can prove the second in the same way.
\end{proof}

{\thm\label{inversemap} It holds that
\[{P_1}^{-1}P_1 =\mathrm{id.}\ \ {on}\ \ \mathcal{R}^{{P_1}^{-1}P_1},\]
\[P_1 {P_1}^{-1}=\mathrm{id.}\ \ {on}\ \ \mathcal{R}^{P_1{P_1}^{-1}}.\]}

\begin{proof}

Let $\pi\in\mathcal{R}^{{P_1}^{-1}P_1}$. Define the sequence\[\lambda_x=\inf_{m\in\mathbb{Z}}\{m\::\:\pi(m)=x\}\ \ \mbox{for}\ x\in\mathbb{Z}\]
with the convention that $\inf \emptyset= \infty$. (See Figure 3,4.) Then, the sequence satisfies one of the following 4 conditions\::

\begin{align*}
&(1)\ \ \cdots<\lambda_{x+1}<\lambda_{x}<\lambda_{x-1}<\cdots     \\
&(2)\ \ -\infty=\lambda_{s}<\lambda_{s-1}<\cdots<\lambda_{x+1}<\lambda_{x}<\lambda_{x-1}<\cdots      \\
&(3)\ \ \cdots<\lambda_{x+1}<\lambda_{x}<\lambda_{x-1}<\cdots<\lambda_{t+1}<\lambda_{t}<\lambda_{t-1}=\infty     \\
&(4)\ \ -\infty=\lambda_{s}<\lambda_{s-1}<\cdots<\lambda_{x+1}<\lambda_{x}<\lambda_{x-1}<\cdots<\lambda_{t}<\lambda_{t-1}=\infty    
\end{align*}

where $s=\liminf_{n\rightarrow-\infty}\pi(n)$ when it is bounded and $t=\liminf_{n\rightarrow\infty}\pi(n)$ when it is bounded. The condition $\inf_{m\leq n}\pi(m)=\pi(n)\ \ i.o.\ as\ \ n\rightarrow\infty$ implies $s\le t$, and if $s=t$, it is the case that $-\infty=\lambda_{s}=\lambda_{t}<\lambda_{t-1}=\infty$.

If (1)\::\:$n=\lambda_{x}$, for some $x$, it holds that
\[P_1\pi(n)=P_1\pi(\lambda_{x})=-\pi(\lambda_{x})+2\inf_{m\leq0}\pi(m)=-\pi(n)+2\inf_{m\leq0}\pi(m)\]
and also it holds that
\[P_1\pi(\lambda_{x})>P_1\pi(\lambda_{x+1})\ \ \mbox{for any}\ -\infty<\lambda_{x+1}<\lambda_{x}<\infty.\]

If (2)\::\:$-\infty<\lambda_{x+1}<n<\lambda_{x}<\infty$ for some $x$, it holds that

\begin{align*}
P_1\pi(n)=\pi(n)-2\pi(\lambda_{x})+2\inf_{m\leq0}\pi(m)&\geq P_1\pi(\lambda_{x})\\
&=-\pi(\lambda_{x})+2\inf_{m\leq0}\pi(m).
\end{align*}

If (3)\::\:$n<\lambda_{s-1}$, it holds that

\begin{align*}
P_1\pi(n)=\pi(n)-2s+2\inf_{m\leq0}\pi(m)&\geq P_1\pi(\lambda_{s-1})-1\\
&=-(s-1)+2\inf_{m\leq0}\pi(m)-1
\end{align*}

If (4)\::\:$n>\lambda_{t}$, it holds that
\[P_1\pi(n)=\pi(n)-2t+2\inf_{m\leq0}\pi(m).\]
and also it holds that
\[P_1\pi(n)=P_1\pi(\lambda_{t})\ \ i.o.\ as\ \ n\rightarrow\infty.\]

From the above discussion, it holds that 

\[\inf_{m\ge n}P_1\pi(m)=\left\{\begin{array}{ll}
-\pi(n)+2\inf_{m\leq0}\pi(m), & \mbox{if  }(1),\\
-\pi(\lambda_{x})+2\inf_{m\leq0}\pi(m), & \mbox{if  }(2),\\
-s+2\inf_{m\leq0}\pi(m), & \mbox{if  }(3),\\
-t+2\inf_{m\leq0}\pi(m), & \mbox{if  }(4).
\end{array}\right.\]

Therefore, if (1), 
\begin{align*}
&{P_1}^{-1}P_1\pi(n)=P_1\pi(n)-2\inf_{m\geq n}P_1\pi(m)+2\inf_{m\geq 0}P_1\pi(m)\\
={}&\left\{-\pi(n)+2\inf_{m\leq0}\pi(m)\right\}-2\left\{-\pi(n)+2\inf_{m\leq0}\pi(m)\right\}+2\inf_{m\geq 0}P_1\pi(m)\\
={}&\pi(n)-2\inf_{m\leq0}\pi(m)+2\inf_{m\geq 0}P_1\pi(m).
\end{align*}

If (2), 
\begin{align*}
&{P_1}^{-1}P_1\pi(n)=P_1\pi(n)-2\inf_{m\geq n}P_1\pi(m)+2\inf_{m\geq 0}P_1\pi(m)\\
={}&\left\{\pi(n)-2\pi(\lambda_{x})+2\inf_{m\leq0}\pi(m)\right\}-2\left\{-\pi(\lambda_{x})+2\inf_{m\leq0}\pi(m)\right\}+2\inf_{m\geq 0}P_1\pi(m)\\
={}&\pi(n)-2\inf_{m\leq0}\pi(m)+2\inf_{m\geq 0}P_1\pi(m).
\end{align*}

If (3), 
\begin{align*}
&{P_1}^{-1}P_1\pi(n)=P_1\pi(n)-2\inf_{m\geq n}P_1\pi(m)+2\inf_{m\geq 0}P_1\pi(m)\\
={}&\left\{\pi(n)-2s+2\inf_{m\leq0}\pi(m)\right\}-2\left\{-s+2\inf_{m\leq0}\pi(m)\right\}+2\inf_{m\geq 0}P_1\pi(m)\\
={}&\pi(n)-2\inf_{m\leq0}\pi(m)+2\inf_{m\geq 0}P_1\pi(m).
\end{align*}

If (4), 
\begin{align*}
&{P_1}^{-1}P_1\pi(n)=P_1\pi(n)-2\inf_{m\geq n}P_1\pi(m)+2\inf_{m\geq 0}P_1\pi(m)\\
={}&\left\{\pi(n)-2t+2\inf_{m\leq0}\pi(m)\right\}-2\left\{-t+2\inf_{m\leq0}\pi(m)\right\}+2\inf_{m\geq 0}P_1\pi(m)\\
={}&\pi(n)-2\inf_{m\leq0}\pi(m)+2\inf_{m\geq 0}P_1\pi(m).
\end{align*}

Therefore it is enough to show that
\[\inf_{m\leq0}\pi(m)=\inf_{m\geq 0}P_1\pi(m),\]
and it is obtained by following inequalities\::
\begin{align*}
\inf_{m\geq 0}P_1\pi(m)&=\inf_{m\geq 0}\left\{\pi(m)-2\inf_{l\leq m}\pi(l)+2\inf_{l\leq 0}\pi(l)\right\}\\
&\geq\inf_{m\geq 0}\left\{\pi(m)-\left(\inf_{l\leq 0}\pi(l)+\inf_{0\leq l\leq m}\pi(l)\right)+2\inf_{l\leq 0}\pi(l)\right\}\\
&=\inf_{m\geq 0}\left\{\pi(m)-\inf_{0\leq l\leq m}\pi(l)\right\}+\inf_{l\leq0}\pi(l)\\
&\geq\inf_{l\leq0}\pi(l).
\end{align*}
On the other hand, by the conditions on $\mathcal{R}^{{P_1}^{-1}P_1}$, there exists $m_1\geq0$ such that $\pi(m_1)=\inf_{l\leq m_1}\pi(l)=\inf_{l\leq 0}\pi(l)$, then
\begin{align*}
\inf_{m\geq 0}P_1\pi(m)&=\inf_{m\geq 0}\left\{\pi(m)-2\inf_{l\leq m}\pi(l)+2\inf_{l\leq 0}\pi(l)\right\}\\
&\leq \pi(m_1)-2\inf_{l\leq m_1}\pi(l)+2\inf_{l\leq 0}\pi(l)\\
&=\inf_{l\leq 0}\pi(l).
\end{align*}

We can prove the second claim in the same way. 
\end{proof}

\begin{figure}
\centering
\scalebox{0.40}{\includegraphics{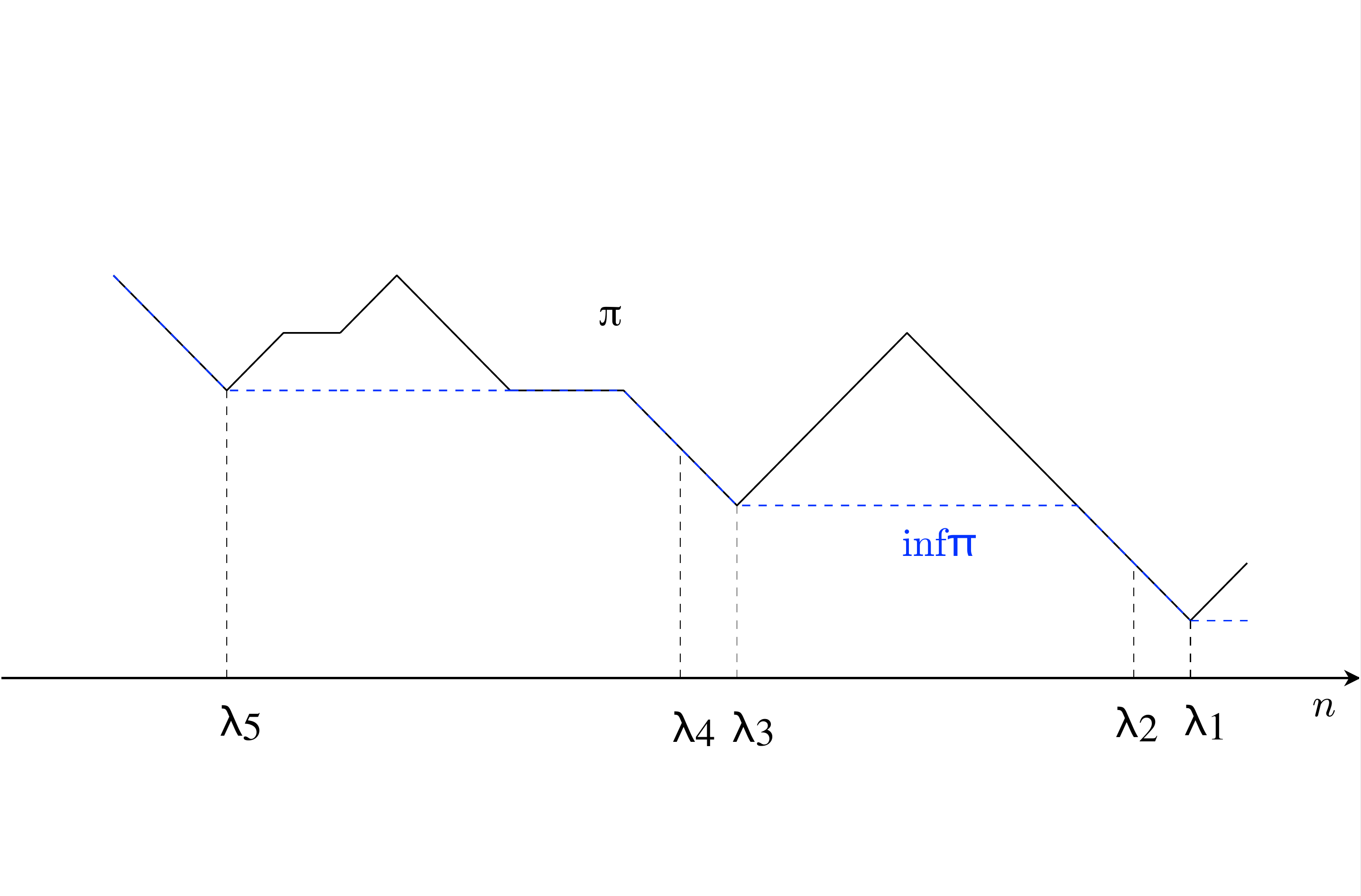}}
\vspace{-20pt}
\caption{Example of the sequence $\{\lambda_x\}$ with $\pi(n),\ \inf_{m\leq n}\pi(m)$.}
\end{figure}

\begin{figure}
\centering
\scalebox{0.40}{\includegraphics{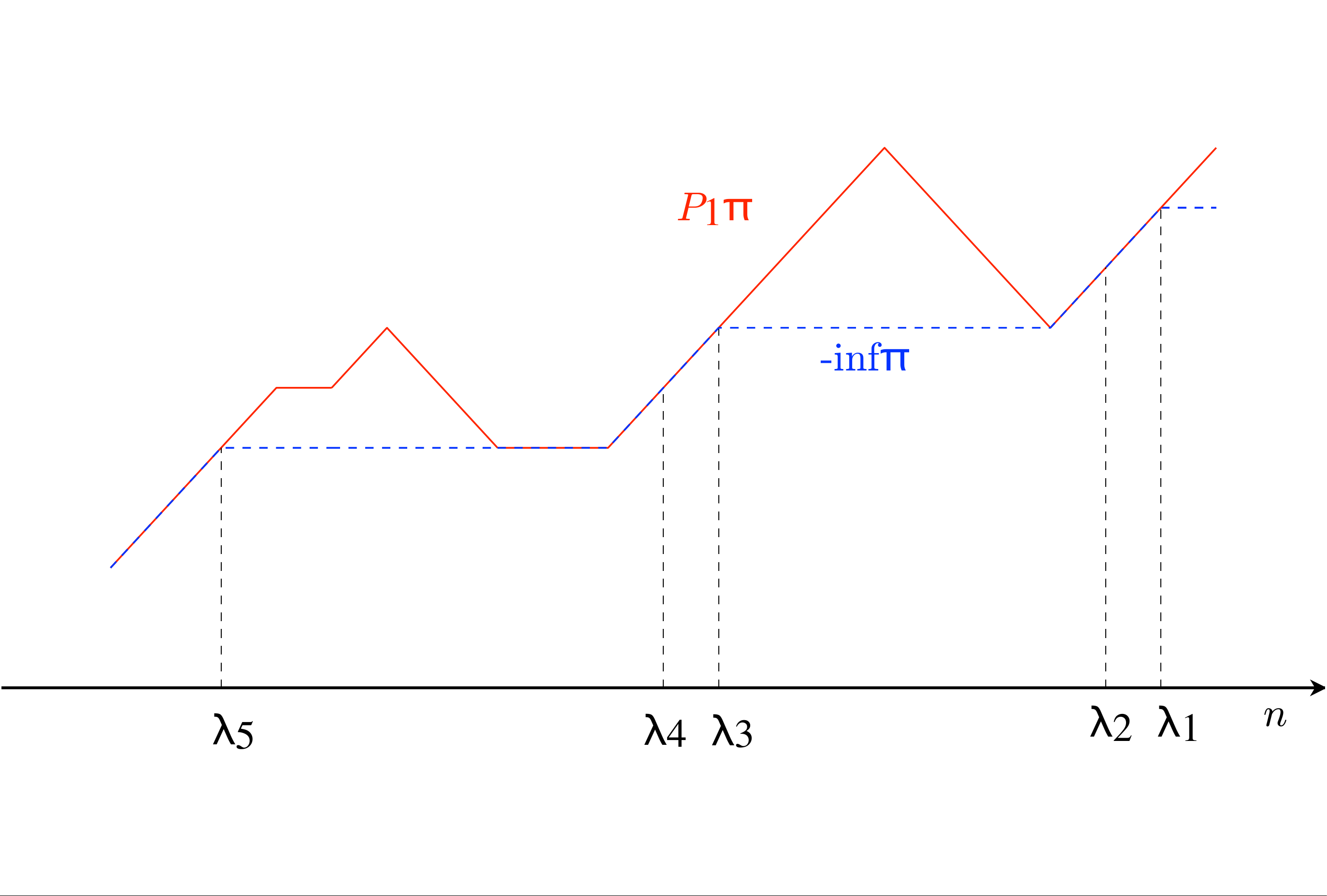}}
\vspace{-20pt}
\caption{The sequence $\{\lambda_x\}$ in figure 3 with $P_1\pi(n)$.}
\end{figure}

{\rem 
The condition $|\pi(n+1)-\pi(n)|\in\{0,1\}$ in $\mathcal{R}^{{P_1}^{-1}P_1}$ and  $\mathcal{R}^{P_1{P_1}^{-1}}$ can be replaced by $|\pi(n+1)-\pi(n)|\in\{0,c\}$ with any positive constant $c$ for Theorem \ref{inversemap} to hold.}

{\rem The condition 
\begin{equation}\label{pi}
\inf_{m\leq n}\pi(m)=\pi(n)\ \ i.o.\ as\ \ n\rightarrow\infty
\end{equation}
in $\mathcal{R}^{{P_1}^{-1}P_1}$ is necessary for ${P_1}^{-1}P_1 =\mathrm{id}$. Indeed, one can check that if $\pi$ does not satisfy \eqref{pi}, the increment of $-\inf_{m\leq n}\pi(m)$ does not match that of $\inf_{m\geq n}P_1\pi(m)$. (See Figure 5, 6.)}

\begin{figure}
\centering
\scalebox{0.40}{\includegraphics{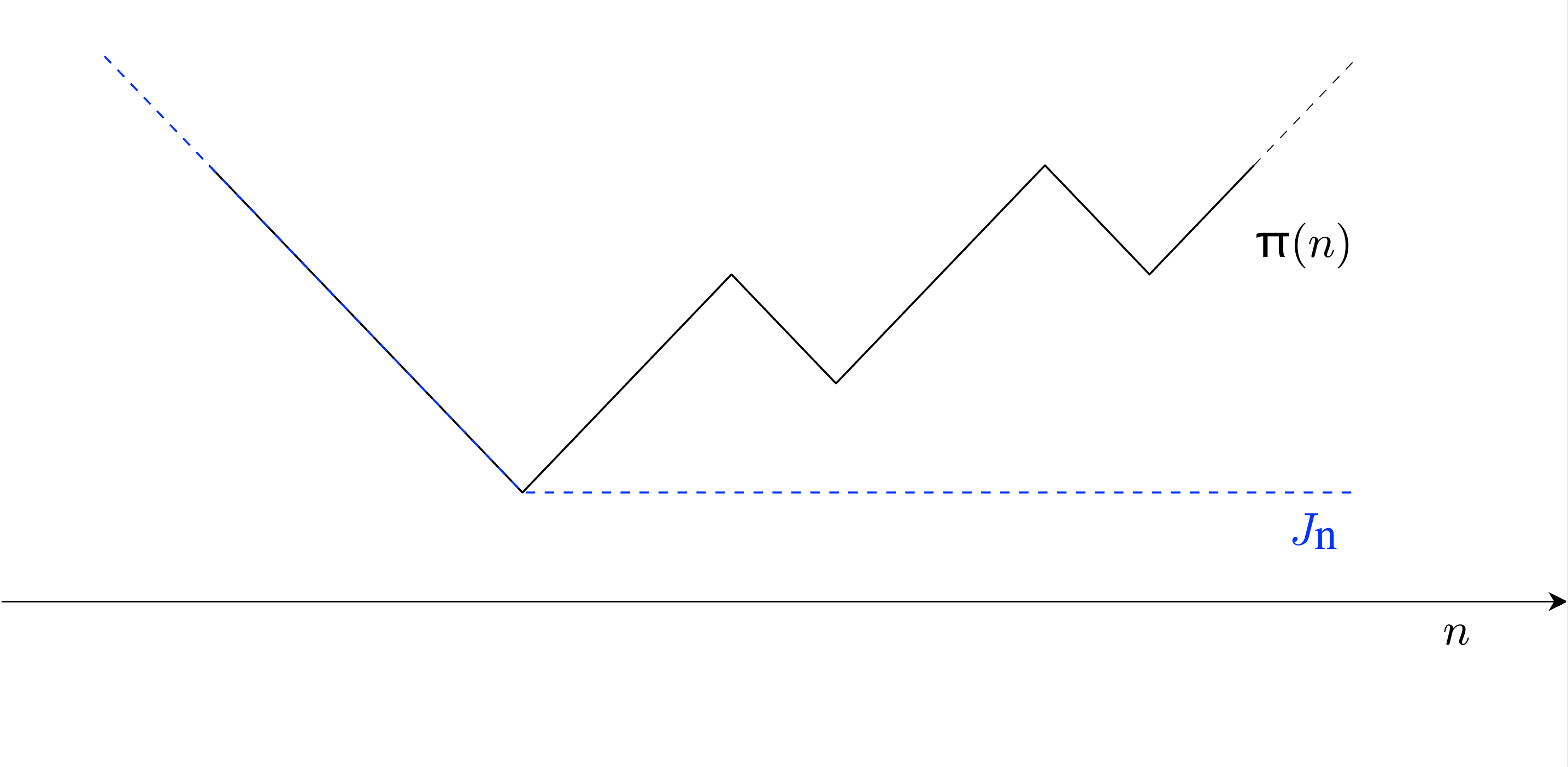}}
\vspace{-30pt}
\caption{Example of $\pi$ not satisfying \eqref{pi} and $J_n:=\inf_{m\leq n}\pi(m)$.}
\end{figure}

\begin{figure}
\centering
\scalebox{0.40}{\includegraphics{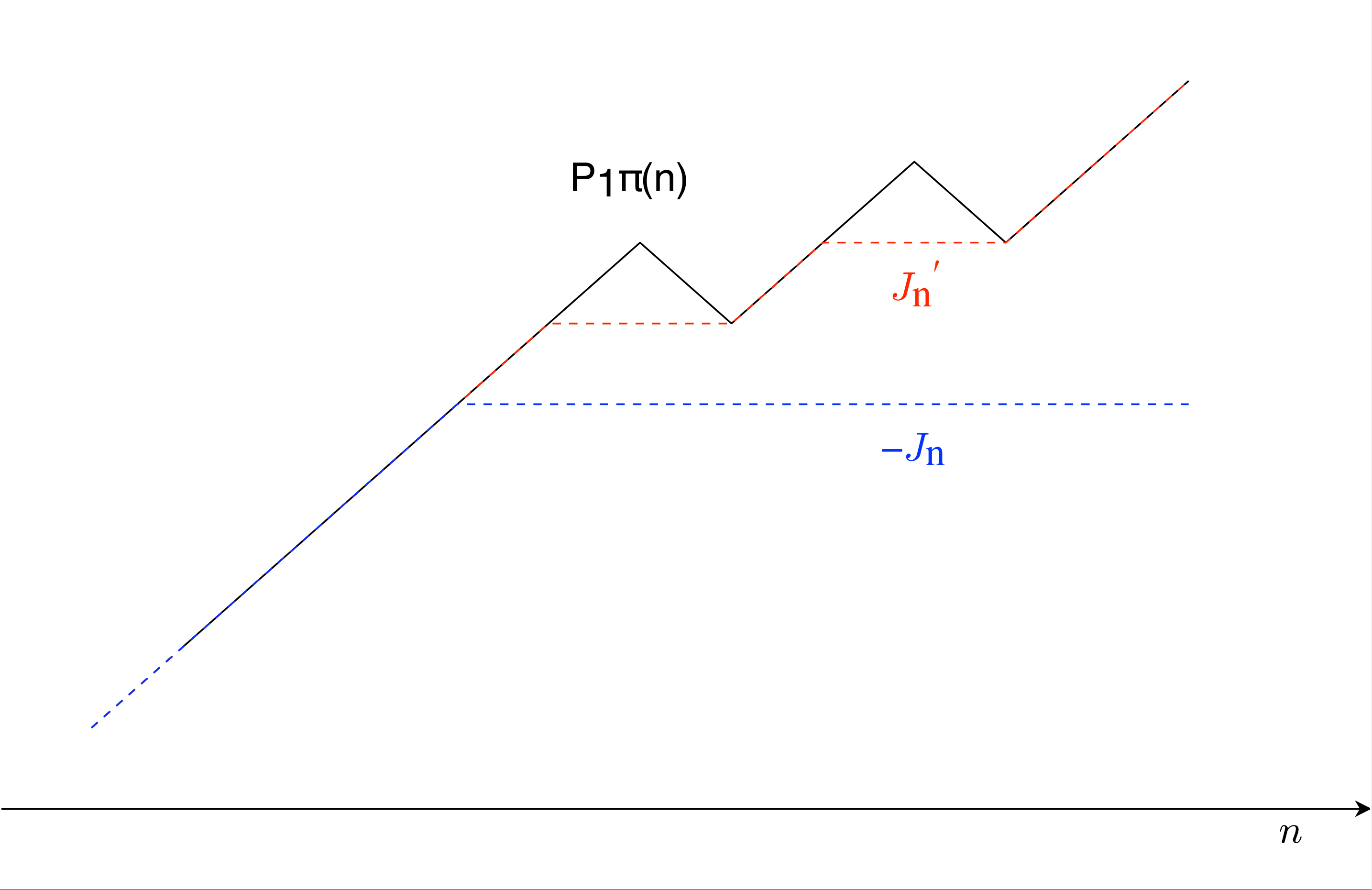}}
\vspace{-10pt}
\caption{Example of $\pi$ not satisfying \eqref{pi} and $J'_n:=\inf_{m\geq n}P_1\pi(m)$.}
\end{figure}

{\cor \label{k-dim inverse} By Remark \ref{orthogonal}}, it holds that
\[{P_\alpha}^{-1}P_\alpha \pi=\pi,\ \ if\ \ \pi_\alpha\in\mathcal{R}^{{P_1}^{-1}P_1},\]
\[P_\alpha {P_\alpha}^{-1}\pi=\pi,\ \ if\ \ \pi_\alpha\in\mathcal{R}^{P_1{P_1}^{-1}},\]
where $\pi_\alpha(n)=\frac{\alpha\cdot\pi(n)}{|\alpha|^2}$.

\section{Path encodings of the multicolor BBS}\label{pathsec}

In the original paper \cite{DS}, the particle configuration is corresponded to the nearest-neighbour walk path $S$ on $\mathbb{Z}$ in $\mathbb{R}$, satisfying $S_0=0$ and $S_n-S_{n-1}=1$ if $\eta_n=0$ and $S_n-S_{n-1}=-1$ if $\eta_n=1$. 
In this section, we extend this concept to the multicolor BBS with $\kappa$-color balls by considering the path $S$ in $\mathbb{R}^{\kappa}$ (Section \ref{one-sided configuration}). In particular, $S$ satisfies $S_0=0$ and $S_n-S_{n-1}=e_i$ if $\eta_n=i\in\{0,1,\cdots,\kappa\}$, where the vectors $e_0,\cdots e_\kappa\in\mathbb{R}^\kappa$ is obtained in Section \ref{Vectors for path encodings}. Then we consider the dynamics of the one-sided multicolor BBS in terms of the ‘carrier’ processes which pick up and drop a certain color ball moving on $\mathbb{Z_+}$ (Section \ref{carrier}), and Pitman transform on $S$ which describes the action $T_i$ (Section \ref{one-sided BBS}). In Section \ref{two-sided BBS}, we extend them to the case of two-sided multicolor BBS. Also we describe the inverse $T^{-1}_i$ and define the reversible set of $S$ for color $i$ such that $T^{-1}_iT_iS=T_iT^{-1}_iS=S$ (Section \ref{inverse}). Moreover, we investigate the set of configurations for which the actions $T_1,T_2,\cdots T_\kappa$ are well-defined and reversible for all times. (Section \ref{invariant set}).

From this section, we fix $\kappa\in\mathbb{N}$ the number of all colors and define the set of numbers representing colors $\mathcal{C}:=\{1,\cdots,\kappa\}$. 

\subsection{Vectors for path encodings}\label{Vectors for path encodings}

In this subsection, we introdece a set of vectors which will be used for path encoding of the particle configuration.

{\df\label{vectorsdef}
Let vectors $e_0,e_1,\cdots,e_\kappa\in\mathbb{R}^{\kappa}$ represent the vertices of a regular $\kappa$-dimensional simplex center the origin, satisfying following conditions\::

\begin{equation}\label{length}
|e_i|=1\ \ \ \forall i\in\mathcal{C}\cup{\{0\}}.
\end{equation}

\begin{equation}\label{product}
e_{i}\cdot e_{j}=-\frac{1}{\kappa}\ \ \ \forall i,j\in\mathcal{C}\cup{\{0\}},i\ne j.
\end{equation}

}

\begin{figure}[H]
\centering
\scalebox{0.40}{\includegraphics{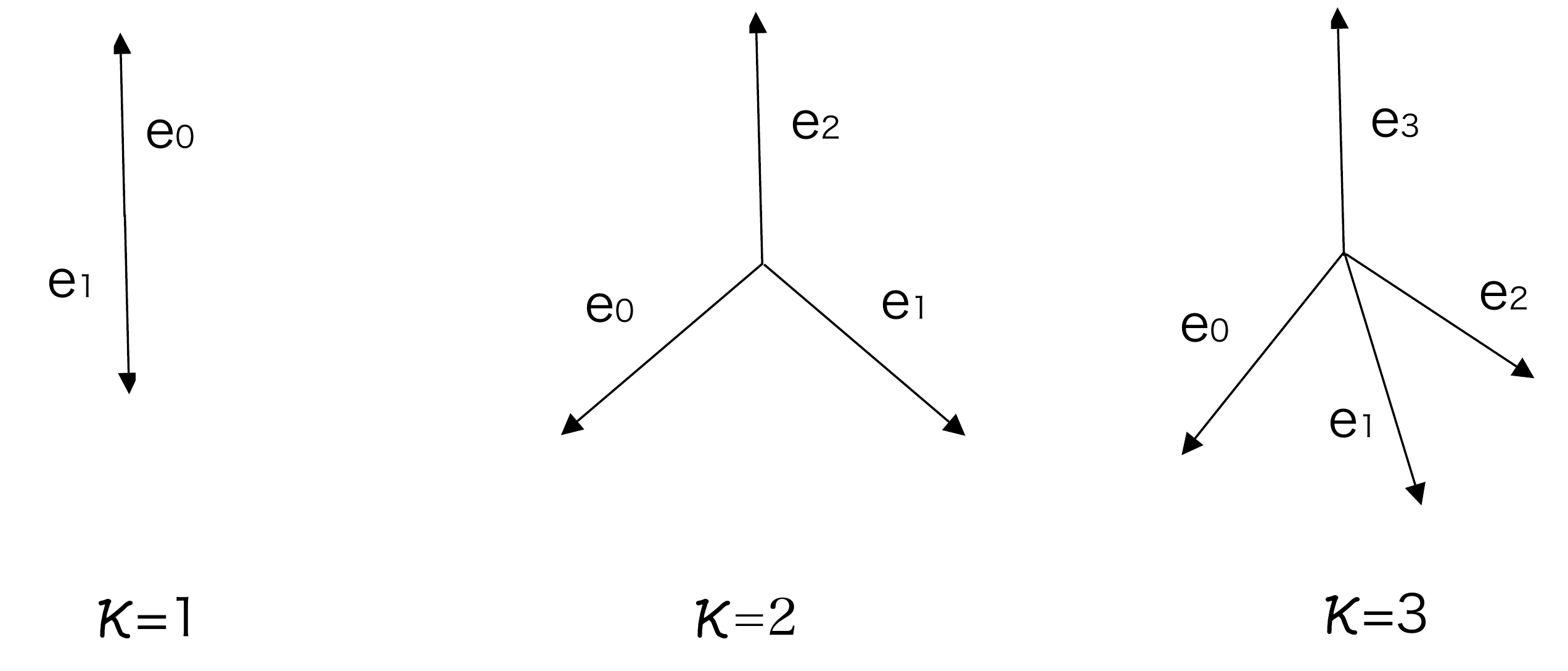}}
\vspace{0pt}
\caption{$e_0,\,e_1\in\mathbb{R}$,\ $e_0,\,e_1,\,e_2\in\mathbb{R}^2$,\ $e_0,\,e_1,\,e_2,\,e_3\in\mathbb{R}^3$}

\end{figure}

{\prop\label{vectorsproperty}
The vectors $e_0,e_1,\cdots,e_\kappa$ have following properties, immediately obtained from \eqref{length} and \eqref{product}, which will be useful in subsequent sections when it comes to defining the path encodings of the particle configuration and considering the actions of the multicolor BBS.

\begin{enumerate}
\setlength{\itemsep}{0.3cm}
\renewcommand{\labelenumi}{(\roman{enumi})}
\item $e_0+e_1+\cdots+e_\kappa=0$
\item Let $a_i\in\mathbb{R}\ $for $i\in\mathcal{C}\cup{\{0\}}$. It holds that
\[a_0e_0+a_1e_1+\cdots+a_\kappa e_\kappa=0 \Leftrightarrow a_0=a_1=\cdots=a_\kappa.\]
\item Let $a_i,a'_i\in\mathbb{R}\ $for $i\in\mathcal{C}\cup{\{0\}}$. Suppose that
\[a_0e_0+a_1e_1+\cdots+a_\kappa e_\kappa=a'_0e_0+a'_1e_1+\cdots+a'_\kappa e_\kappa.\]
Then there is a constant $c$ such that $a_i=a'_i+c$ for any $i$. In addition, suppose that
\[a_0+a_1+\cdots a_\kappa=a'_0+a'_1+\cdots +a'_\kappa.\]
Then it is the case that $a_i=a'_i$ for any $i$.
\item Let $a_l\in\mathbb{R}$ for $l\in\mathcal{C}\cup{\{0\}}$, and $d_j\in\mathbb{R}$ for $j\in\mathcal{C}$. It holds that\begin{align*}
&a_0e_0+a_1e_1+\cdots+a_\kappa e_\kappa=d_i(e_i-e_0)+\sum_{j\in\mathcal{C},j\ne i}d_je_j\\
\Leftrightarrow{}&d_j=a_j-\frac{a_0+a_i}{2}\ \ \forall j\in\mathcal{C} 
\end{align*}
for any $i\in\mathcal{C}$.
\item Any set of $\kappa$ vectors in $\{e_0,e_1,\cdots,e_\kappa\}$ is the basis of $\mathbb{R}^\kappa$.
\item For any $v\in\mathbb{R}^\kappa$, there is an $\kappa+1$-tuple $a_0,\cdots,a_\kappa$ of real numbers satisfying
\[v=a_0e_0+\cdots+a_\kappa e_\kappa,\ \ a_0+\cdots+a_\kappa=0\]
\end{enumerate}

}

\subsection{Configuration of the one-sided multicolor BBS}\label{one-sided configuration}

In this section, we consider the one-sided multicolor BBS, and denote the particle configuration by $\eta=(\eta_n)_{n\in\mathbb{\kappa}}\in \{0,1,2,\cdots,\kappa\}^{\mathbb{N}}$ As in the introduction, we write $\eta_n=i$ if there is a particle of color $i\in\mathcal{C}$ at site $n$, and $\eta_n=0$ if there is no particle at site $n$. 

We define a nearest-neighbour path in $\mathbb{R}^{\kappa}$ as the path encoding of a particle configuration. 

{\df\label{half path encoding} Given the particle configuration by $\eta=(\eta_n)_{n\in\mathbb{\kappa}}\in \{0,1,2,\cdots,\kappa\}^{\mathbb{N}}$, we define $S=(S_n)_{n\in\mathbb{Z}_+}$ by setting

\begin{equation}\label{increment}
S_0=0\ \ S_n-S_{n-1}=e_i\ \ \text{if}\ \ \eta_n=i.
\end{equation}
The S is called the path encoding of $\eta$. We can describe it as 
\begin{equation}\label{path encoding}
S_n=a_0(n)e_0+a_1(n)e_1+\cdots+a_\kappa(n)e_\kappa\ \ 
\end{equation}
for $n\in\mathbb{Z_{+}}$, where $a_i(n)\in\mathbb{Z_+},\ i\in\mathcal{C}$ is the number of the particles of color $i$ at the sites located from $1$ to $n$, $a_0(n)\in\mathbb{Z_+}$ is the number of the empty sites located from $1$ to $n$, and $a_i(0)=0,\ i\in\mathcal{C}\cup{\{0\}}$. Also we define the path space in $\mathbb{R}^\kappa$ as follows\::
\[\mathcal{S}_+:=\{S:\mathbb{Z_+}\rightarrow\mathbb{R}^{\kappa}\::S_0=0,\ S_{n+1}-S_{n}\in\{e_0,e_1,\cdots,e_\kappa\},\ \forall n\in\mathbb{Z_+}\}.\]
}

{\exmp For $\eta=(0,1,1,2,\cdots)$, the path encoding S is given by\[S_0=0,\ S_1=e_0,\ S_2=e_0+e_1,\ S_3=e_0+2e_1,\ S_4=e_0+2e_1+e_2,\ \cdots\]}

{\rem By the definition, it clearly holds that the map from $\eta=(\eta_n)_{n\in\mathbb{N}}\in \{0,1,2,\cdots,\kappa\}^{\mathbb{N}}$ to $S\in\mathcal{S}_+$ is one to one. Also it holds that\[a_0(n)+a_1(n)+\cdots+a_\kappa(n)=n\ \ \forall n.\]Therefore, from Proposition \ref{vectorsproperty} $($\hspace{-1pt}ⅲ\hspace{-1pt}$)$, the map from $(a_0(n),a_1(n),\cdots,a_\kappa(n))\in\mathbb{Z_+}^{\kappa+1}$ to $S_n\in\mathbb{R}^\kappa$ is one to one for any $n\in\mathbb{Z_+}$.}

\vspace{10pt}

For the subsequent sections, we introduce some operators of $S$.


{\df \label{Adef1}
For $i\in\mathcal{C}$, we define the function $A_i:\mathcal{S}_+\rightarrow\mathbb{Z}^{\mathbb{Z_+}}$ given by
\begin{equation}\label{A1}
A_iS_n=a_0(n)-a_i(n), 
\end{equation}
for $S_n=a_0(n)e_0+a_1(n)e_1+\cdots+a_\kappa(n)e_\kappa,\ n\in\mathbb{Z_+}$.
}

{\rem\label{Adef2}
For $S_n=a_0(n)e_0+a_1(n)e_1+\cdots+a_\kappa(n) e_\kappa$, Proposition \ref{vectorsproperty} $($\hspace{-1pt}ⅳ\hspace{-1pt}$)$ shows
\[S_n=\frac{1}{2}\left\{a_i(n)-a_0(n)\right\}(e_i-e_0)+\sum_{j\ne0,i}d_j(n)e_j\]
and \eqref{product} implies
\[e_j\cdot(e_i-e_0)=0\ \ \ \forall i,j\in\mathcal{C},\ i\ne j.\]
Therefore, the projection of $S_n$ along $(e_i-e_0)$ is equal to 
\[\frac{(e_i-e_0)\cdot S_n}{|e_i-e_0|^2}(e_i-e_0)=\frac{1}{2}\left\{a_i(n)-a_0(n)\right\}(e_i-e_0)=-\frac12A_iS_n(e_i-e_0),\]
then it holds that
\begin{equation}\label{A2}
A_iS_n=-2\frac{(e_i-e_0)\cdot S_n}{|e_i-e_0|^2}.
\end{equation}
}

{\rem\label{projection}
From Remark \ref{Adef2}, 
we can write $S_n$ as the sum of the vector projection on $(e_i-e_0)$ and the vector orthogonal to $(e_i-e_0)$ as following
\[S_n=-\frac12A_iS_n(e_i-e_0)+\left(S_n+\frac12A_iS_n(e_i-e_0)\right).\]
Then, by Proposition \ref{half orthogonal}, it holds that 
\begin{align*}
P_{e_i-e_0}S_n&=P_{e_i-e_0}\left(-\frac{1}{2}A_iS_n\left(e_i-e_0\right)+\left(S_n+\frac12A_iS_n(e_i-e_0)\right)\right)\\
&=P_1\left(-\frac{1}{2}A_iS_n\right)\left(e_i-e_0\right)+\left(S_n+\frac12A_iS_n(e_i-e_0)\right).
\end{align*}
}

{\df \label{taudef1}

We define the the permutation operator $\tau_{(0,i)}:\mathcal{S}_+\rightarrow \mathcal{S}_+$ given by
\begin{equation}\label{tau1}
\tau_{(0,i)}S_n=a_i(n)e_0+a_0(n)e_i+\sum_{j\ne0,i}a_j(n)e_j
\end{equation}
for $S_n=a_0(n)e_0+a_1(n)e_1+\cdots+a_\kappa(n)e_\kappa,\ n\in\mathbb{Z_+}$.
}

{\rem \label{taudef2}

Comparing
\begin{align*}
S_n&=\frac{1}{2}\left\{a_i(n)-a_0(n)\right\}(e_i-e_0)+\sum_{j\ne0,i}d_j(n)e_j\\
&=-\frac12A_iS_n(e_i-e_0)+\left(S_n+\frac12A_iS_n(e_i-e_0)\right),
\end{align*}
and
\begin{align*}
\tau_{(0,i)}S_n&=\frac{1}{2}\left\{a_0(n)-a_i(n)\right\}(e_i-e_0)+\sum_{j\ne0,i}d_j(n)e_j\\
&=\frac12A_iS_n(e_i-e_0)+\left(S_n+\frac12A_iS_n(e_i-e_0)\right),
\end{align*}
it is the case that $\tau_{(0,i)}$ is the operator which multiply only the vector projection part of $S$ along $e_i-e_0$ by $-1$. Also it holds that 
\begin{equation}\label{tau2}
\tau_{(0,i)}S_n=S_n+A_iS_n(e_i-e_0).
\end{equation}

}


\subsection{Carrier process for the one-sided multicolor BBS}\label{carrier}

We introduce the concept of carrier with respect to particles of a certain color $i\in\mathcal{C}$. It moves along $\mathbb{Z_+}$ from left to right picking up a particle of color $i$ when it crosses one, and dropping off a particle of color $i$ when it is holding at least one particle and sees an empty site. The dynamic $T_i$ can be viewed in terms of this carrier. The carrier process is given as follows.

{\df The carrier process $W^{(i)}=\{W^{(i)}_n\}_{n\in\mathbb{Z_+}}$ of the color $i$ associated with $\eta\in \{0,1,2,\cdots,\kappa\}^{\mathbb{N}}$ is defined by $W^{(i)}_0=0$ and

\begin{equation}\label{W}
W^{(i)}_n=\left\{\begin{array}{ll}
W^{(i)}_{n-1}+1, & \mbox{if }\eta_n=i,\\
W^{(i)}_{n-1}, & \mbox{if }\eta_n=j,\ j\ne 0,i\\
W^{(i)}_{n-1}, & \mbox{if }\eta_n=0\mbox{ and }W^{(i)}_{n-1}=0,\\
W^{(i)}_{n-1}-1, & \mbox{if }\eta_n=0\mbox{ and }W^{(i)}_{n-1}>0.
\end{array}\right.
\end{equation}}

$W$ is obtained from $S$ as following lemma.

{\lem\label{W=M-S} It holds that\[W^{(i)}_n=\sup_{0\leq m\leq n}A_iS_m-A_iS_n,\ \ \forall n\in\mathbb{Z_+}.\]}

\begin{proof}

We prove it by induction. Clearly the result is true for $n=0$. Suppose that $W^{(i)}_{n-1}=\sup_{0\leq m\leq n-1}A_iS_m-A_iS_{n-1}$ for some $n\geq1$. 

Now, if $\eta_n=i$, then $A_iS_n=A_iS_{n-1}-1$ and $\sup_{0\leq m\leq n}A_iS_m=\sup_{0\leq m\leq n-1}A_iS_m$, and so
\[\left\{\sup_{0\leq m\leq n}A_iS_m-A_iS_n\right\}-\left\{\sup_{0\leq m\leq n-1}A_iS_m-A_iS_{n-1}\right\}=1.\]
If $\eta_n=j,\ j\ne 0,i$, then $A_iS_n=A_iS_{n-1}$ and $\sup_{0\leq m\leq n-1}A_iS_m=\sup_{0\leq m\leq n}A_iS_m$, and so
\[\left\{\sup_{0\leq m\leq n}A_iS_m-A_iS_n\right\}-\left\{\sup_{0\leq m\leq n-1}A_iS_m-A_iS_{n-1}\right\}=0.\]
Moreover, if $\eta_n=0$ and $W^{(i)}_{n-1}=0$, then it is the case that $\sup_{0\leq m\leq n}A_iS_m=A_iS_n$, and so 
\[\left\{\sup_{0\leq m\leq n}A_iS_m-A_iS_n\right\}-\left\{\sup_{0\leq m\leq n-1}A_iS_m-A_iS_{n-1}\right\}=0.\]
Similarly, if $\eta_n=0$ and $W^{(i)}_{n-1}>0$, then $A_iS_n=A_iS_{n-1}+1$ and $\sup_{0\leq m\leq n}A_iS_m=\sup_{0\leq m\leq n-1}A_iS_m$, and so 
\[\left\{\sup_{0\leq m\leq n}A_iS_m-A_iS_n\right\}-\left\{\sup_{0\leq m\leq n-1}A_iS_m-A_iS_{n-1}\right\}=-1.\]
Thus it holds that 
\[W^{(i)}_{n}-W^{(i)}_{n-1}=\left\{\sup_{0\leq m\leq n}A_iS_m-A_iS_n\right\}-\left\{\sup_{0\leq m\leq n-1}A_iS_m-A_iS_{n-1}\right\}\]
which by the inductive hypothesis implies 
\[W^{(i)}_{n}=\sup_{0\leq m\leq n}A_iS_m-A_iS_n.\]

\end{proof}

\subsection{Action of the carrier for the one-sided multicolor BBS}\label{one-sided BBS}



In this section, we consider the action $T_i$ on S given by \eqref{path encoding}. We fix the color $i\in\mathcal{C}$. 
From the viewpoint of the carrier process, we can write $T_i$ as 
\[T_i(\eta)_n=\mathbf{1}_{\{W^{(i)}_n=W^{(i)}_{n-1}-1\}},\ \ \forall n\in\mathbb{N}.\]
For $j\ne i$, the numbers $\{a_j(n)\}_{n\in\mathbb{Z_{+}}}$ do not change under the action $T_i$, so the path encoding $T_iS=(T_iS_n)_{n}\in\mathbb{Z_{+}}$ of $T_i\eta$ can be described as follows,
\[T_iS_n=a'_0(n)e_0+a'_i(n)e_i+\sum_{j\ne0,i}a_j(n)e_j\]
for some $a'_0(n)$ and $a'_i(n)$.

Then $T_i$ satisfies the following formula.

{\lem \label{2M-S}It holds that \[a'_0(n)-a'_i(n)=2\sup_{0\leq m\leq n}\left\{a_0(m)-a_i(m)\right\}-\left\{a_0(n)-a_i(n)\right\}.\]
That is, from Definition \ref{Adef1}
\begin{align*}
A_iT_iS_n&=2\sup_{0\leq m\leq n}A_iS_m-A_iS_n\\
&=P_1(-A_iS)_n
\end{align*}
by using Pitman transform (\ref{P_1}) in Definition \ref{1dim oneside pitman}.}

\begin{proof}

It is easy to check that 
\[2\mathbf{1}_{\{A_iS_n-A_iS_{n-1}=1\}}=1+(A_iS_n-A_iS_{n-1})-\mathbf{1}_{\{\eta_n\ne0,i\}}.\]
This equation and Theorem \ref{W=M-S} show that
\begin{align*}
& A_iT_iS_n-A_iT_iS_{n-1}\\
={}&1-2\mathbf{1}_{\{W^{(i)}_n=W^{(i)}_{n-1}-1\}}-\mathbf{1}_{\{\eta_n\ne0,i\}}\\
={}&1-2\mathbf{1}_{\{A_iS_{n-1}<\sup_{0\leq m\leq n-1}A_iS_m,\ A_iS_n-A_iS_{n-1}=1\}}-\mathbf{1}_{\{\eta_n\ne0,i\}}\\
={}&1-\left(2\mathbf{1}_{\{A_iS_n-A_iS_{n-1}=1\}}-2\mathbf{1}_{\{A_iS_{n-1}=\sup_{0\leq m\leq n-1}A_iS_m,\ A_iS_n-A_iS_{n-1}=1\}}\right)-\mathbf{1}_{\{\eta_n\ne0,i\}}\\
={}&-(A_iS_n-A_iS_{n-1})+2\mathbf{1}_{\{A_iS_{n-1}=\sup_{0\leq m\leq n-1}A_iS_m,\ A_iS_n-A_iS_{n-1}=1\}}.
\end{align*}
Summing over the increments, we obtain 

\begin{align*}
&A_iT_iS_n-A_iT_iS_0\\
={}&\sum^n_{m=1}\left(A_iT_iS_m-A_iT_iS_{m-1}\right)\\
={}&A_iS_0-A_iS_n+2\sum^n_{m=1}\mathbf{1}_{\{A_iS_{n-1}=\sup_{0\leq m\leq n-1}A_iS_m,\ A_iS_n-A_iS_{n-1}=1\}}\\
={}&A_iS_0-A_iS_n+2\left(\sup_{0\leq m\leq n}A_iS_m-\sup_{0\leq m\leq 0}A_iS_m\right).
\end{align*}
Since $A_iS_0=A_iT_iS_0=\sup_{0\leq m\leq 0}A_iS_m=0$, the claim is proved.

\end{proof}

{\thm \label{oneside TP}It holds that \[T_iS=\tau_{(0,i)}P_{e_i-e_0}S,\ \forall S\in\mathcal{S_+}\]
where P is one-sided Pitman transform defined in Definition \ref{oneside pitman}.}

\begin{proof}
By Remark \ref{projection} and Lemma \ref{2M-S}, it holds that

\begin{align*}
P_{e_i-e_0}S_n&=P_1\left(-\frac{1}{2}A_iS_n\right)\left(e_i-e_0\right)+\sum_{j\ne0,i}d_j(n)e_j\\
&=\frac{1}{2}P_1\left(-A_iS_n\right)\left(e_i-e_0\right)+\sum_{j\ne0,i}d_j(n)e_j\\
&=\frac{1}{2}A_iT_iS_n\left(e_i-e_0\right)+\sum_{j\ne0,i}d_j(n)e_j
\end{align*}
where $d_j(n)=a_j(n)-\frac{a_0(n)+a_i(n)}{2}=a_j(n)-\frac{a'_0(n)+a'_i(n)}{2},\ j\ne0,i$.

On the other hand, by Definition \ref{taudef1},
\begin{align*}
\tau_{(0,i)}T_iS_n&=a'_i(n)e_0+a'_0(n)e_i+\sum_{j\ne0,i}a_j(n)e_j\\
&=\frac{1}{2}(a'_0(n)-a'_i(n))\left(e_i-e_0\right)+\sum_{j\ne0,i}d_j(n)e_j\\
&=\frac{1}{2}A_iT_iS_n\left\{e_i-e_0\right\}+\sum_{j\ne0,i}d_j(n)e_j.
\end{align*}
Therefore, we obtain the equation $\tau_{(0,i)}T_iS_n=P_{e_i-e_0}S_n$ for any $n\in\mathbb{Z_+}$.
\end{proof}

{\rem The dynamic T for the 1-color BBS in the paper \cite{DS} is expressed as follows\::
\[TS_n=2\sup_{0\leq m\leq n}S_m-S_n,\]
where $S_n=a_0(n)e_0+a_1(n)e_1=a_0(n)-a_1(n)$. This is also called Pitman transform and corresponds to Lemma \ref{2M-S}. For the multicolor case, however, the supremum expression is
\begin{align*}
&2\sup_{0\leq m\leq n}\frac{(e_0-e_i)\cdot S_m}{|e_0-e_i|^2}(e_0-e_i)-S_n\\
={}&2\sup_{0\leq m\leq n}\frac{(e_0-e_i)\cdot S_m}{|e_0-e_i|^2}(e_0-e_i)-\left(a_0(n)e_0+a_1(n)e_1+\cdots+a_\kappa(n)e_\kappa\right)\\
={}&2\sup_{0\leq m\leq n}\frac{a_0(m)-a_i(m)}{2}(e_0-e_i)-\left(\frac{a_0(n)-a_i(n)}{2}(e_0-e_i)+\sum_{j\ne0,i}d_j(n)e_j\right)\\
={}&\frac12\left(2\sup_{0\leq m\leq n}\left\{a_0(m)-a_i(m)\right\}-\left\{a_0(n)-a_i(n)\right\}\right)(e_0-e_i)-\sum_{j\ne0,i}d_j(n)e_j\\
={}&\frac{a'_0(m)-a'_i(m)}{2}(e_0-e_i)-\sum_{j\ne0,i}d_j(n)e_j, 
\end{align*}
where $d_j(n)=a_j(n)-\frac{a_0(n)+a_i(n)}{2}=a_j(n)-\frac{a'_0(n)+a'_i(n)}{2},\ j\ne0,i$. Then this does not correspond to $T_iS_n$ because the sign of $d_j(n)$ is negative. 
This is the reason why we use infimum expression of Pitman transform.

}

{\rem From Theorem \ref{oneside TP}, it holds that

\begin{align*}
T_2T_1&=\left(\tau_{(0,2)}P_{e_2-e_0}\right)\left(\tau_{(0,1)}P_{e_1-e_0}\right)\\
&=\tau_{(0,1)}\left(\tau_{(1,2)}P_{e_2-e_1}\right)P_{e_1-e_0}\\
&=\tau_{(0,1)}\tau_{(1,2)}P_{e_2-e_1}P_{e_1-e_0}.
\end{align*}
Similarly, the dynamic $T$ of the multicolor BBS is as follows\::
\[T=T_\kappa\cdots T_2T_1=\tau_{(0,1)}\tau_{(1,2)}\cdots\tau_{(\kappa-1,\kappa)}P_{e_{\kappa}-e_{\kappa-1}}\cdots P_{e_2-e_1}P_{e_1-e_0}.\]
}

\subsection{Two-sided multicolor BBS}\label{two-sided BBS}

In this section, we extend the particle configuration to $\eta=(\eta_n)_{n\in\mathbb{Z}}\in \{0,1,2,\cdots,\kappa\}^{\mathbb{Z}}$.

We can again obtain the path encoding $S=(S_n)_{n\in\mathbb{Z}}$ of the $\eta$ given by (\ref{increment}) and (\ref{path encoding}). In this case, for $i\in\mathcal{C}$ and $n\geq1$, $a_i(n)$ means the number of the particles of color $i$ at the sites located from $1$ to $n$, and, for $i\in\mathcal{C}$ and $n\leq-1$, $-a_i(n)$ means the same at the sites located from $n+1$ to $0$. The same is true for the number of the empty sites. Also we define that $a_i(0)=0$ for $i\in\mathcal{C}\cup{\{0\}}$. As in the case of one-sided multicolor BBS, it obviously holds $a_0(n)+a_1(n)+\cdots+a_\kappa(n)=n\ \ \forall n\in\mathbb{Z}$. 

Also we define the path space in $\mathbb{R}^\kappa$\::
\[\mathcal{S}^0:=\{S=(S_n)_{n\in\mathbb{Z}}\::\:S_0=0,\ S_{n+1}-S_{n}\in\{e_0,e_1,\cdots,e_\kappa\},\ \forall n\in\mathbb{Z}\}.\]
Moreover, we define the function $A_i:\mathcal{S}^0\rightarrow\mathbb{R}^{\mathbb{Z}}$ and the operator $\tau_{(0,i)}:\mathcal{S}^0\rightarrow \mathcal{S}^0$ given by \eqref{A1} and \eqref{tau1}.

Whilst in the one-sided case, carrier process W and the actions $T_i,\ i=1,\cdots \kappa$ are defined for any $S\in\mathcal{S_+}$ (that is, for any configuration $\eta\in\{0,1,2,\cdots,\kappa\}^{\mathbb{N}}$), in the two-sided case, the following restriction on $S$ is required to define the carrier and actions\::
\begin{equation}\label{twoside condition}
\limsup_{n\rightarrow-\infty}A_iS_n<\infty.
\end{equation}
This condition can be transformed as follows\::
\begin{align*}
\limsup_{n\rightarrow-\infty}A_iS_n<\infty&\Leftrightarrow \sup_{n\leq0}A_iS_n<\infty\\
&\Leftrightarrow \sup_{n\leq0}\left\{a_0(n)-a_i(n)\right\}<\infty\\
&\Leftrightarrow \inf_{n\leq0}\left\{\left(-a_0(n)\right)-\left(-a_i(n)\right)\right\}>-\infty\\
&\Leftrightarrow -A_iS\in\mathcal{R}^{P_1}
\end{align*}
and this means that the number of particles of color $i$ is not too much compared with the number of empty sites in the left side.

Indeed, in section 2.4 in the paper \cite{DS}, two-sided multicolor BBS is understood with two-sided carrier process 
\[W^{(i)}_n=\sup_{m\leq n}A_iS_m-A_iS_n\]
under the condition (\ref{twoside condition}).

Also the path encoding $T_iS_n=a'_0(n)e_0+a'_i(n)e_i+\sum_{j\ne0,i}a_j(n)e_j$ of $T_i\eta$, is obtained by the equation
\begin{equation}\label{twoside action}
A_iT_iS_n=2\sup_{m\leq n}A_iS_m-A_iS_n-2\sup_{m\leq 0}A_iS_m
\end{equation}
under the condition (\ref{twoside condition}). Then, in the same way as proof of Theorem \ref{oneside TP}, it holds that
\[T_iS=\tau_{(0,i)}P_{e_i-e_0}S\]
where P is two-sided Pitman transform defined in Definition \ref{twoside pitman}.

From the above discussion, the next set is obtained\::
\begin{align*}
\mathcal{S}^{T_i}:&=\{S\in\mathcal{S}^0\::\:T_iS\mbox{ well-defined}\}\\
&=\{S\in\mathcal{S}^0\::\:\limsup_{n\rightarrow-\infty}A_iS_n<\infty\}
\end{align*}

\subsection{Inverse of the action}\label{inverse}

In the previous section, we found that $T_i=\tau_{(0,i)}P_{e_i-e_0}$ on $\mathcal{S}^{T_i}$. Then we can defined $T^{-1}_i=P^{-1}_{e_i-e_0}\tau_{(0,i)}$ on an appropriate set, where $P^{-1}$ is defined by Definition \ref{inverse pitman}. 

As in the proof of Theorem \ref{oneside TP}, $T_i$ acts on $S_n$ as follows\::
\begin{align*}
T_iS_n&=\tau_{(0,i)}P_{e_i-e_0}S_n\\
&=\tau_{(0,i)}\left(\frac{1}{2}P_{1}\left(-A_iS\right)_n(e_i-e_0)+\sum_{j\ne0,i}d_j(n)e_j\right).
\end{align*}

Therefore, Theorem \ref{inversemap} shows 

\begin{align*}
T^{-1}_iT_iS_n&=T^{-1}_i\left(\frac{1}{2}P_{1}\left(-A_iS\right)_n(e_i-e_0)+\sum_{j\ne0,i}d_j(n)e_j\right)\\
&=P^{-1}_{e_i-e_0}\tau_{(0,i)}\tau_{(0,i)}\left(\frac{1}{2}P_{1}\left(-A_iS\right)_n(e_i-e_0)+\sum_{j\ne0,i}d_j(n)e_j\right)\\
&=P^{-1}_{e_i-e_0}\left(\frac{1}{2}P_{1}\left(-A_iS\right)_n(e_i-e_0)+\sum_{j\ne0,i}d_j(n)e_j\right)\\
&=P^{-1}_1\left(\frac{1}{2}P_{1}\left(-A_iS\right)\right)_n(e_i-e_0)+\sum_{j\ne0,i}d_j(n)e_j\\
&=\frac{1}{2}P^{-1}_1P_{1}\left(-A_iS\right)_n(e_i-e_0)+\sum_{j\ne0,i}d_j(n)e_j\\
&=\frac{1}{2}\left(-A_iS\right)_n(e_i-e_0)+\sum_{j\ne0,i}d_j(n)e_j\\
&=S_n
\end{align*}
if and only if $-A_iS\in\mathcal{R}^{{P_1}^{-1}P_1}$.

On the other hand, by Remark \ref{taudef2}, $T^{-1}_i$ acts on $S_n$ as follows\::
\begin{align*}
T^{-1}_iS_n&=P^{-1}_{e_i-e_0}\tau_{(0,i)}S_n\\
&=P^{-1}_{e_i-e_0}\tau_{(0,i)}\left(-\frac12\left(A_iS_n\right)(e_i-e_0)+\sum_{j\ne0,i}d_j(n)e_j\right)\\
&=P^{-1}_{e_i-e_0}\left(\frac12\left(A_iS_n\right)(e_i-e_0)+\sum_{j\ne0,i}d_j(n)e_j\right)\\
&=P^{-1}_{1}\left(\frac12A_iS\right)_n(e_i-e_0)+\sum_{j\ne0,i}d_j(n)e_j\\
&=\frac12P^{-1}_{1}\left(A_iS\right)_n(e_i-e_0)+\sum_{j\ne0,i}d_j(n)e_j
\end{align*}

Therefore, Theorem \ref{inverse} and Remark \ref{taudef2} shows

\begin{align*}
T_iT^{-1}_iS_n&=\tau_{(0,i)}P_{e_i-e_0}\left(\frac12P^{-1}_{1}\left(A_iS\right)_n(e_i-e_0)+\sum_{j\ne0,i}d_j(n)e_j\right)\\
&=\tau_{(0,i)}\left(P_1\left(\frac12P^{-1}_{1}\left(A_iS\right)\right)_n(e_i-e_0)+\sum_{j\ne0,i}d_j(n)e_j\right)\\
&=\tau_{(0,i)}\left(\frac12P_1P^{-1}_{1}\left(A_iS\right)_n(e_i-e_0)+\sum_{j\ne0,i}d_j(n)e_j\right)\\
&=\tau_{(0,i)}\left(\frac12A_iS_n(e_i-e_0)+\sum_{j\ne0,i}d_j(n)e_j\right)\\
&=-\frac12A_iS_n(e_i-e_0)+\sum_{j\ne0,i}d_j(n)e_j\\
&=S_n
\end{align*}
if and only if $A_iS\in\mathcal{R}^{P_1{P_1}^{-1}}$.

Above discussion gives the following theorem characterizing the following sets\::

\begin{align*}
\mathcal{S}^{T^{-1}_iT_i}:=&\{S\in\mathcal{S}^0\::\:T_iS,T^{-1}_iTS\mbox{ well-defined},\ T^{-1}_iT_iS=S\}\\
\mathcal{S}^{T_iT^{-1}_i}:=&\{S\in\mathcal{S}^0\::\:T^{-1}_iS,T_iT^{-1}_iS\mbox{ well-defined},\ T_iT^{-1}_iS=S\}.
\end{align*}

{\thm It holds that 
\begin{align*}
\mathcal{S}^{T^{-1}_iT_i}&=\{S\in\mathcal{S}^0\::\:-A_iS\in\mathcal{R}^{{P_1}^{-1}P_1}\}\\
&=\{S\in\mathcal{S}^0\::\:\inf_{m\leq0}\left(-A_iS_m\right)>-\infty,\ \inf_{m\leq n}\left(-A_iS_m\right)=-A_iS_n,\ i.o.\ as\ n\rightarrow\infty\}\\
&=\{S\in\mathcal{S}^0\::\:\sup_{m\leq0}A_iS_m<\infty,\ \sup_{m\leq n}A_iS_m=A_iS_n,\ i.o.\ as\ n\rightarrow\infty\},
\end{align*}
and 
\begin{align*}
\mathcal{S}^{T_iT^{-1}_i}&=\{S\in\mathcal{S}^0\::\:A_iS\in\mathcal{R}^{P_1{P_1}^{-1}}\}\\
&=\{S\in\mathcal{S}^0\::\:\inf_{m\geq0}A_iS_m>-\infty,\ \inf_{m\geq n}A_iS_m=A_iS_n,\ i.o.\ as\ n\rightarrow-\infty\}.
\end{align*}
}

{\rem The above conditions can be transformed as follows\::
\[\sup_{m\leq n}A_iS_m=A_iS_n,\ i.o.\ as\ n\rightarrow\infty\ \Leftrightarrow\ \sup_{n\in\mathbb{Z}}A_iS_n=\limsup_{n\rightarrow\infty}A_iS_n,\]
\[\inf_{m\geq n}A_iS_m=A_iS_n,\ i.o.\ as\ n\rightarrow-\infty\ \Leftrightarrow\ \inf_{n\in\mathbb{Z}}A_iS_n=\liminf_{n\rightarrow-\infty}A_iS_n.\]

Then it holds that

\[\mathcal{S}^{T^{-1}_iT_i}=\{S\in\mathcal{S}^0\::\:M^{(i)}_0<\infty,\ \limsup_{n\rightarrow\infty}A_iS_n=M^{(i)}_\infty\},\]
\[\mathcal{S}^{T_iT^{-1}_i}=\{S\in\mathcal{S}^0\::\:I^{(i)}_0>-\infty,\ \liminf_{n\rightarrow-\infty}A_iS_n=I^{(i)}_{-\infty}\}.\]
where, we define
\[M^{(i)}_0:=\sup_{n\leq 0}A_{(i)}S_n,\ M^{(i)}_\infty:=\sup_{n\in\mathbb{Z}}A_{(i)}S_n,\]
\[I^{(i)}_0:=\inf_{n\geq 0}A_iS_n,\ I^{(i)}_{-\infty}:=\inf_{n\in\mathbb{Z}}A_iS_n.\]

Also we obtain the following set\::

\begin{align*}
\mathcal{S}^{rev}_i:&=\{S\in\mathcal{S}^0\::\:T_iS,T^{-1}_iS,T^{-1}_iTS,T_iT^{-1}_iS\mbox{ well-defined},\ T^{-1}_iT_iS=T_iT^{-1}_iS=S\}\\
&=\{S\in\mathcal{S}^0\::\:M^{(i)}_0<\infty,\ I^{(i)}_0>-\infty,\ \limsup_{n\rightarrow\infty}A_iS_n=M^{(i)}_\infty,\ \liminf_{n\rightarrow-\infty}A_iS_n=I^{(i)}_{-\infty}\}.
\end{align*}
}

\subsection{Set of configurations}\label{invariant set}

Even if $S\in\mathcal{S}^{rev}_i$ holds, it does not necessarily hold that $T_iS\in\mathcal{S}^{rev}_i$. In the paper \cite{DS} for the 1-color BBS, the set

\[\mathcal{S}^{inv}_i:=\{S\in\mathcal{S}^0\::\:T^{k}_iS\in\mathcal{S}^{rev}_i,\ \forall k\in\mathbb{Z}\}\]
is characterized as following lemma. 

{\lem \label{Sinv}
For any $i\in\mathcal{C}$, it holds that 
\[\mathcal{S}^{inv}_i=\bigcup_{*_1,*_2\in\{sub-critical(i),critical(i)\}}\left(\mathcal{S}_{*_1}^-\cap\mathcal{S}_{*_2}^+\right),\]
where
\begin{align*}
\mathcal{S}_{sub-critical(i)}^{\pm}&:=\left\{S\in\mathcal{S}^0\::\:\lim_{n\rightarrow\pm\infty}\frac{A_iS_n}{F_i(n)}=1,\ \exists F_i\in\mathcal{F}\right\},\\
\mathcal{S}_{critical(i)}^{\pm}&:=\left\{S\in\mathcal{S}^0\::\: \sup_{n\in\mathbb{Z}}W^{(i)}_n<\infty,\:\limsup_{n \to \pm \infty}A_iS_n= \liminf_{n \to \pm \infty}A_iS_n + \sup_{n}W^{(i)}_n \in \mathbb{R}\right\},\\
\mathcal{F}&:=\{F:\mathbb{Z}\rightarrow\mathbb{R}\::\:\mbox{increasing function, }\lim_{n\rightarrow\infty} F(n)=\infty,\ \lim_{n\rightarrow-\infty} F(n)=-\infty\}.
\end{align*}
Moreover, it holds that
\begin{equation}\label{Fnochange}
\lim_{n\rightarrow\pm\infty}\frac{A_iS_n}{F_i(n)}=1,\ \exists F_i\in\mathcal{F}\ \Rightarrow\ \lim_{n\rightarrow\pm\infty}\frac{A_iT_iS_n}{F_i(n)}=1
\end{equation}
and 

\begin{equation}\label{FisM}
\lim_{n\rightarrow\pm\infty}\frac{A_iS_n}{F_i(n)}=1,\ \exists F_i\in\mathcal{F}\ \Leftrightarrow\ \lim_{n\rightarrow\pm\infty}\frac{A_iS_n}{\sup_{m\leq n}A_iS_m}=1.
\end{equation}
}

For the study of the multicolor BBS theory, it is natural to ask when $T^{-1}TS=TT^{-1}S=S$ is true where $T$ is any composition of $T_i,\ i\in\mathcal{C}$ such as $T=T_\kappa\cdots T_2T_1,\ T=T_2T_1T^2_2$ etc. In other words, what is the condition for $S$ to be in the following set ? 

\[\mathcal{S}^{inv}_{\mathcal{C}}:=\{S\in\mathcal{S}^0\::\:TS\in\bigcap_{i\in\mathcal{C}}\mathcal{S}^{rev}_i\mbox{\ for any composition }T\mbox{ of }T_i,\ i\in\mathcal{C}\}.\]
One might expect that 
\[\mathcal{S}^{inv}_{\mathcal{C}}\supseteq\bigcap_{i\in\mathcal{C}}\mathcal{S}^{inv}_i\]
but this is not true. (See Remark \ref{noinvariant}.) 
The main result of this section is the following theorem which gives a sufficient condition for $S$ to be in the set $S\in\mathcal{S}^{inv}_{\mathcal{C}}$.

{\thm \label{invariantthm}
Define the subset of $\bigcap_{i\in\mathcal{C}}\left(\mathcal{S}_{sub-critical(i)}^-\cap\mathcal{S}_{sub-critical(i)}^+\right)$ such that $F_i$ and $F_j$ have the same asymptotic behavior as $n\rightarrow\pm\infty$ for any $i,j\in\mathcal{C}$ as follows, 
\[\mathcal{S}^{good}_{\mathcal{C}}:=\left\{S\in\mathcal{S}^0\::\:\forall i\in\mathcal{C}\ \exists F_i\in\mathcal{F},\ \lim_{n\rightarrow\pm\infty}\frac{A_iS_n}{F_i(n)}=1\mbox{ and }\limsup_{n\rightarrow\pm\infty}\frac{F_j(n)}{F_i(n)}<\infty\ \forall i,j\in\mathcal{C}\right\}.
\]
It holds that
\[\mathcal{S}^{inv}_{\mathcal{C}}\supseteq\mathcal{S}^{good}_{\mathcal{C}}.\]
}

To prove the above result, we prepare a simple lemma.

{\lem For any $i, j\in\mathcal{C},\ i\ne j$, and $S\in\mathcal{S}^{T_i}$it holds that

\begin{equation}\label{a'1}
A_jT_iS_n=A_jS_n+ W^{(i)}_n-M^{(i)}_0
\end{equation}
and
\begin{equation}\label{a'2}
A_jT_iS_n=A_jS_n+\frac12\left(A_iT_iS_n-A_iS_n\right)
\end{equation}
for any $n\in\mathbb{Z}$}.

\begin{proof}

Let $S_n=a_0(n)e_0+a_1(n)e_1+\cdots+a_\kappa(n)e_\kappa$ and $T_iS_n=a'_0(n)e_0+a'_i(n)e_i+\sum_{k\ne0,i}a_k(n)e_k$. Then (\ref{twoside action}) shows
\[a'_0(n)-a'_i(n)=2\sup_{m\leq n}\left\{a_0(m)-a_i(m)\right\}-\left\{a_0(n)-a_i(n)\right\}-2M^{(i)}_0\]
By adding $a'_0(n)+a'_i(n)=a_0(n)+a_i(n)$ to the above equation, we have
\[2a'_0(n)=2\sup_{m\leq n}\left\{a_0(m)-a_i(m)\right\}+2a_i(n)-2M^{(i)}_0.\]
Then it follows that
\[a'_0(n)=a_0(n)+\sup_{m\leq n}\left\{a_0(m)-a_i(m)\right\}-\left\{a_0(n)-a_i(n)\right\}-M^{(i)}_0\]
Since $A_jT_iS_n=a'_0(n)-a_j(n)$ and $\sup_{m\leq n}A^{(i)}S_m-A^{(i)}S_n=W^{(i)}_n$, the first claim is proved. Also $a'_0(n)+a'_i(n)=a_0(n)+a_i(n)$ shows 
\[2a'_0(n)-\left\{a'_0(n)-a'_i(n)\right\}=2a_0-\left\{a_0(n)-a_i(n)\right\}\]
then,
\[a'_0(n)=a_0(n)+\frac12\left(A_iT_iS_n-A_iS_n\right)\]
and this prove the second claim.
\end{proof}

\begin{proof}[Proof of Theorem \ref{invariantthm}]

Suppose that $S\in\mathcal{S}^{good}_{\mathcal{C}}$. It is enough to show that $T_iS\in\mathcal{S}^{good}_{\mathcal{C}}$ for any $i\in\mathcal{C}$, so for that we show

\[\lim_{n\rightarrow\pm\infty}\frac{A_jT_iS_n}{F_j(n)}=1\]
for any $i,j\in\mathcal{C}$. From \eqref{a'2}, we can write
\[\frac{A_jT_iS_n}{F_j(n)}=\frac{A_jS_n}{F_j(n)}+\frac12\frac{F_j(n)}{F_i(n)}\left(\frac{A_iT_iS_n}{F_i(n)}-\frac{A_iS_n}{F_i(n)}\right).\]
By the assumption and \eqref{Fnochange},  it holds that
\[\lim_{n\rightarrow\pm\infty}\frac{A_jS_n}{F_j(n)}=1,\ \lim_{n\rightarrow\pm\infty}\frac{A_iS_n}{F_i(n)}=1,\ \lim_{n\rightarrow\pm\infty}\frac{A_iT_iS_n}{F_i(n)}=1.\]
Then the condition $\limsup_{n\rightarrow\pm\infty}\frac{F_j(n)}{F_i(n)}<\infty$ shows the conclusion.

\end{proof}

{\rem\label{noinvariant}

Now we consider three examples of the configurations with $\mathcal{C}=\{1,2\}$. Each example shows one of the following three claims.

\begin{align*}
(a)\ \ \mathcal{S}^{inv}_{\mathcal{C}}&\not\supseteq \bigcap_{i\in\mathcal{C}}\left(\mathcal{S}_{sub-critical(i)}^-\cap\mathcal{S}_{sub-critical(i)}^+\right),\\
(b)\ \ \mathcal{S}^{inv}_{\mathcal{C}}&\not\supseteq \bigcap_{i\in\mathcal{C}}\left(\mathcal{S}_{critical(i)}^-\cap\mathcal{S}_{critical(i)}^+\right),\\
(c)\ \ \mathcal{S}^{inv}_{\mathcal{C}}&\varsupsetneq\mathcal{S}^{good}_{\mathcal{C}}.
\end{align*}

\vspace{10pt}

$(a)$\ We give an example of $\eta$ whose path encoding $S$ satisfies 
\[S\in\bigcap_{i\in\mathcal{C}}\left(\mathcal{S}_{sub-critical(i)}^-\cap\mathcal{S}_{sub-critical(i)}^+\right),\ T_2S\notin\mathcal{S}_{sub-critical(1)}^+,\ T_2S\notin\mathcal{S}_{critical(1)}^+\]
Let $\eta$ be as follows\::

\vspace{5pt}
\makebox[0.37cm][r]{$\eta=\ $}$(\cdots\ 0\ \eta_0=0\ 0\ 2_{(1)}\ (0\ 1)_{(1)}\ 0\ (0\ 1)_{(2)}\ 0\ 2_{(3)}\ (0\ 1)_{(3)}\ 0\ (0\ 1)_{(4)}\ 0\ 2_{(5)}\ (0\ 1)_{(5)}\ 0\ (0\ 1)_{(6)}\cdots$

\vspace{5pt}
\makebox[2.6cm][r]{}$\cdots\cdots\ 0\ 2_{(2m-1)}\ (0\ 1)_{(2m-1)}\ 0\ (0\ 1)_{(2m)}\ \cdots),$
\vspace{5pt}

\noindent where $i_{(k)}:=i\ i\ \cdots i$ means k consecutive i, and $(i\ j)_k:=i\ j\ i\ j\ \cdots\ i\ j$ means that i and j alternately appear k times. For simplicity, Figure 6 and Figure 7 show the graph of $A_2S_n$ and $A_1S_n$ skipping places where there is no increase or decrease where $S$ is path encoding of $\eta$. As seen in Figure 6, it holds that 

\begin{align*}
1&\geq\limsup_{n\rightarrow\infty}\frac{A_2S_n}{\sup_{m\leq n}A_2S_m}\\
&\geq\liminf_{n\rightarrow\infty}\frac{A_2S_n}{\sup_{m\leq n}A_2S_m}\\
&=\lim_{k\rightarrow\infty}\frac{1+(1+3)+(1+5)+\cdots+(1+2k-1)-(2k-1)}{1+(1+3)+(1+5)+\cdots+(1+2k-1)}\\
&=1.
\end{align*}
As seen in Figure 7, it holds that $\left|\sup_{m\leq n}A_1S_m-A_1S_n\right|\leq1,\ \forall n$ and $\lim_{n\rightarrow\infty}A_1S_n=\infty$, then $\lim_{n\rightarrow\infty}\frac{A_1S_n}{\sup_{m\leq n}A_1S_m}=1$. Also it clearly holds that $\lim_{n\rightarrow-\infty}\frac{A_1S_n}{\sup_{m\leq n}A_1S_m}=\lim_{n\rightarrow-\infty}\frac{A_2S_n}{\sup_{m\leq n}A_2S_m}=1$. Therefore, by Lemma \ref{Sinv}, $S\in\bigcap_{i\in\mathcal{C}}\left(\mathcal{S}_{sub-critical(i)}^-\cap\mathcal{S}_{sub-critical(i)}^+\right)$.

However, the configuration of $T_2\eta$ is as follows\::

\vspace{5pt}
\makebox[0.37cm][r]{$T_2\eta=\ $}$(\cdots\ 0\ \eta_0=0\ 0\ 0_{(1)}\ (2\ 1)_{(1)}\ 0\ (0\ 1)_{(2)}\ 0\ 0_{(3)}\ (2\ 1)_{(3)}\ 0\ (0\ 1)_{(4)}\ 0\ 0_{(5)}\ (2\ 1)_{(5)}\ 0\ (0\ 1)_{(6)}\cdots$

\vspace{5pt}
\makebox[2.6cm][r]{}$\cdots\cdots\ 0\ 0_{(2m-1)}\ (2\ 1)_{(2m-1)}\ 0\ (0\ 1)_{(2m)}\ \cdots).$

\vspace{5pt}
\noindent Figure 8, the graph of $A_1T_2S_n$, shows that $\liminf_{n\rightarrow\infty}\frac{A_1T_2S_n}{\sup_{m\leq n}A_1T_2S_m}=\frac12$. Therefore, $T_2S\notin\mathcal{S}_{sub-critical(1)}^+$. Also $T_2S\notin\mathcal{S}_{critical(1)}^+$ is obvious. Then, from Lemma \ref{Sinv}, $T_2S\notin\mathcal{S}^{inv}_1$.
Such a phenomenon occurs because $W^{(2)}_n$ can be arbitrarily large and it causes a gap between the asymptotic behavior of $A_1T_2S_n$ and that of $A_1S_n$ as $n\rightarrow\infty$ from the equation
$A_1T_2S_n=A_1S_n+ W^{(2)}_n-M^{(2)}_0$ by \eqref{a'1}. 

\vspace{5pt}

\begin{figure}[h]
\centering
\scalebox{0.30}{\includegraphics{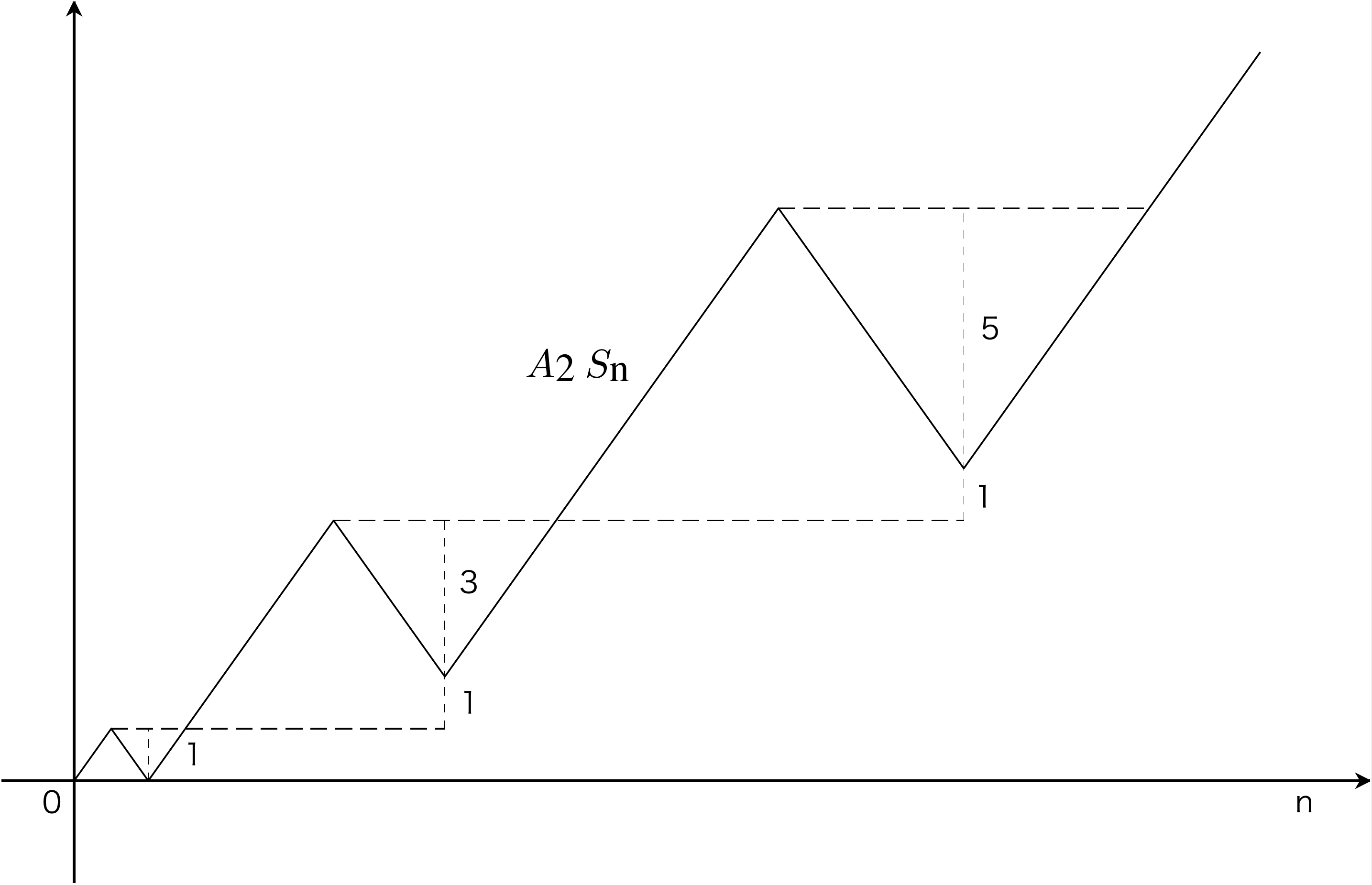}}
\vspace{4pt}
\caption{The graph of $A_2S_n$ skipping places where there is no increase or decrease.}
\end{figure}

\vspace{6pt}
\begin{figure}[h]
\centering
\scalebox{0.30}{\includegraphics{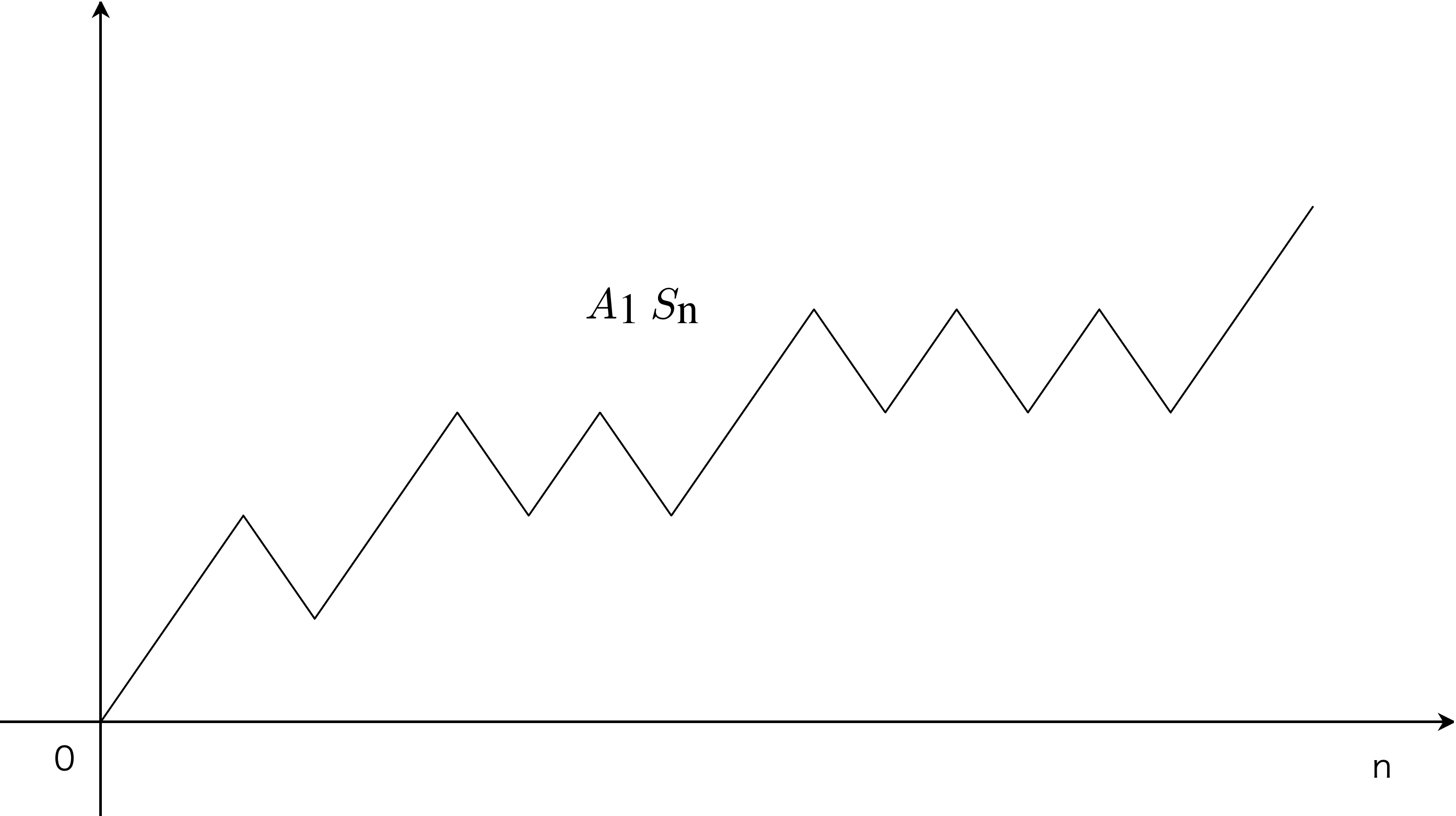}}
\caption{The graph of $A_1S_n$ skipping places where there is no increase or decrease.}
\end{figure}

\vspace{6pt}
\begin{figure}[h]
\centering
\scalebox{0.30}{\includegraphics{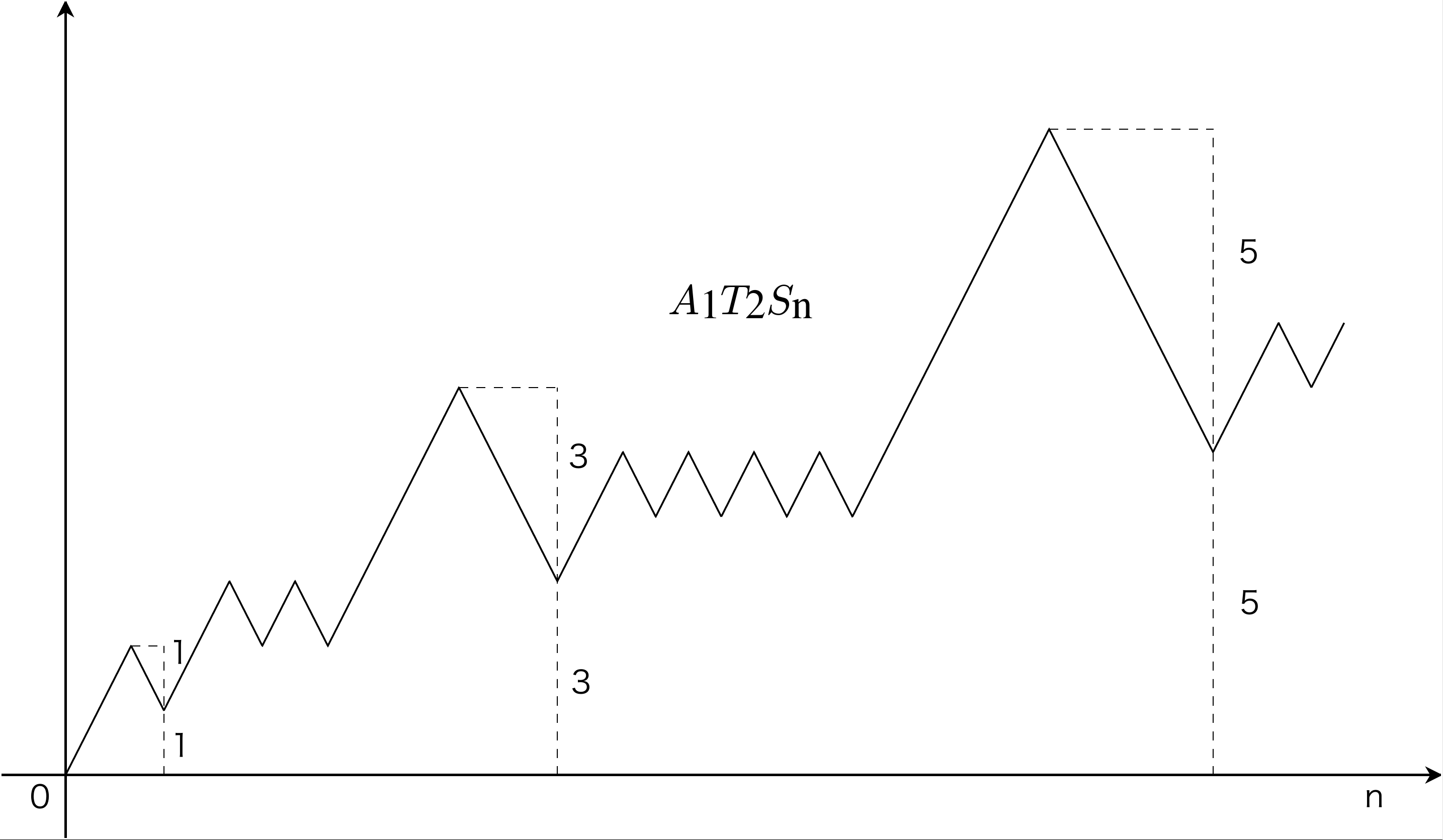}}
\vspace{4pt}
\caption{The graph of $A_1T_2S_n$ skipping places where there is no increase or decrease.}
\end{figure}

$(b)$\ We give an example of $\xi$ whose path encoding $S^{(\xi)}$ satisfies 
\[S^{(\xi)}\in\bigcap_{i\in\mathcal{C}}\left(\mathcal{S}_{critical(i)}^-\cap\mathcal{S}_{critical(i)}^+\right),\ T_2S^{(\xi)}\notin\mathcal{S}_{sub-critical(1)}^+,\ T_2S^{(\xi)}\notin\mathcal{S}_{critical(1)}^+.\]
Let $\xi$ be as follows\::

\vspace{5pt}
\makebox[1.85cm][r]{$\xi=\ $}$(\cdots\ 0\ 1\ 2\ 0\ 1\ 2\ 0\ 2\ 1\ 0\ 1\ 2\ 0\ 1\ 2\ 0\ 1\ 2\ \cdots).$

\vspace{5pt}

\noindent Then,

\vspace{5pt}
\makebox[1.85cm][r]{$T_2\xi=\ $}$(\cdots\ 2\ 1\ 0\ 2\ 1\ 0\ 2\ 0\ 1\ 2\ 1\ 0\ 2\ 1\ 0\ 2\ 1\ 0\ \cdots)$

\vspace{5pt}

\noindent and they show above conditions.

$(c)$\ We give an example of $\zeta$ whose path encoding $S^{(\zeta)}$ satisfies 
\[S^{(\zeta)}\in\bigcap_{i\in\mathcal{C}}\left(\mathcal{S}_{critical(i)}^-\cap\mathcal{S}_{critical(i)}^+\right),\ S^{(\zeta)}\in\mathcal{S}^{inv}_\mathcal{C}.\]
Let $\zeta$ be as follows\::

\vspace{5pt}
\makebox[1.85cm][r]{$\zeta=\ $}$(\cdots\ 0\ 1\ 2\ 0\ 1\ 2\ 0\ 1\ 2\ 0\ 1\ 2\ 0\ 1\ 2\ 0\ 1\ 2\ \cdots).$

\vspace{5pt}

\noindent$TS^{(\zeta)}\in\bigcap_{i\in\mathcal{C}}\left(\mathcal{S}_{critical(i)}^-\cap\mathcal{S}_{critical(i)}^+\right)$, where $T$ is any composition of $T_1$ and $T_2$, because the configuration of $T\zeta$ is always repeating $(012)$ or $(021)$. Therefore, it holds that $S^{(\zeta)}\in\mathcal{S}^{inv}_{\mathcal{C}}$.

}

\section{Random initial configurations}
In this section, we consider the case when the initial configuration is random. 
Suppose that $\eta=(\eta_n)_{n\in\mathbb{Z}}$ is an ergodic sequence which is stationary with respect to the space shift. In particular, if we assume that the densities of the balls of color $i$  

\begin{equation}\label{density}
p_i=\mathbf{P}(\eta_0=i)<p_0=\mathbf{P}(\eta_0=0),\ \ \ \forall i\in\mathcal{C},
\end{equation}
then ergodicity implies that $A_iS$ satisfies

\[\frac{A_iS}{n}=\frac{a_0(n)-a_i(n)}{n}\rightarrow p_0-p_i>0,\ \ \ \mathbf{P}\mbox{-a.s.}\]
as $n\rightarrow\pm\infty$. Thus we obtain the following result, which yields that $(T^kS)_{k\in\mathbb{Z}}$ is well-defined and reversible by Theorem \ref{invariantthm}.

{\lem \label{ergodic}If $\eta=(\eta_n)_{n\in\mathbb{Z}}$ is a stationary, ergodic sequence satisfying \eqref{density}, then it holds that\[\frac{A_iS}{(p_0-p_i)n}\rightarrow 1,\ \ \ \mathbf{P}\mbox{-a.s.}\]
as $n\rightarrow\pm\infty$ for any $i\in\mathcal{C}$. In particular, $S\in\mathcal{S}^{good}_C,\ \mathbf{P}\mbox{-a.s.}$.}

Next, it is natural for random initial configuration to ask whether the law of $\eta$ is preserved by $T_i$, that is, $T_i\eta\buildrel{d}\over{=}\eta$. We introduce the example of an invariant measure in Section \ref{iid}. Moreover, we consider generalized multicolor BBS whose dynamic is defined for continuous path in $\mathbb{R}^\kappa$, and generalize each object appearing in the discrete case (Section \ref{onR}). And in Section \ref{BMD}, we check that $\kappa$-dimensional Brownian motion with certain drift is invariant under the action of the multicolor BBS, and it is obtained by appropriate scaling limit of asymmetric random walk with distribution $\mathbf{P}(S_m-S_{m-1}=e_j)=\frac1{\kappa+1}+\frac{c_j}{\sqrt{n\kappa}},j\in\mathcal{C}\cup{0}$ such that $c_0>c_i,\forall i\in\mathcal{C}$, which represents a high density particle configuration.

\subsection{Independent and identically distributed initial configuration}\label{iid}

Suppose that $\eta=(\eta_n)_{n\in\mathbb{Z}}$ is given by a sequence of i.i.d. random variables with following distribution
\begin{equation}\label{iidprob}
p_i=\mathbf{P}(\eta_0=i)<p_0=\mathbf{P}(\eta_0=0),\ \ \ \forall i\in\mathcal{C},
\end{equation}
then it satisfies \eqref{density} and the conditions in Lemma \ref{ergodic}. Furthermore, $S$ is a random walk path in $\mathbb{R}^{\kappa}$ satisfying $S_0=0$ and 
\[\mathbf{P}(S_n-S_{n-1}=e_j)=p_j,\ \ \ \forall j\in\mathcal{C}\cup\{0\},\]
where the increments of S are independent. 

{\thm \label{iidinvariant} If $\eta=(\eta_n)_{n\in\mathbb{Z}}$ is given by a sequence of i.i.d. random variables with \eqref{iidprob}, it holds that \[T_i\eta\buildrel{d}\over{=}\eta\]
for any $i\in\mathcal{C}$.}

\begin{proof}

We introduce some notations.

For each $n\in\mathbb{Z}$, define transform $f_n:\{0,1,\cdots,\kappa\}^{\mathbb{Z}}\rightarrow\{1,\cdots,\kappa\}$ given by 
\[f_n(\eta)=\left\{\begin{array}{ll}
\eta_n, & \mbox{if  }\eta_n\not\in\{0,i\},\\
i, & \mbox{if  }\eta_n\in\{0,i\}.\\
\end{array}\right.\]
Then it holds that $\left(f_n(T_i\eta)\right)_{n\in\mathbb{Z}}=\left(f_n(\eta)\right)_{n\in\mathbb{Z}}$, for each $\eta\in\{0,1,\cdots,\kappa\}^{\mathbb{Z}}$.

For each $n\in\mathbb{Z}$ and $\eta\in\{0,1,\cdots,\kappa\}^{\mathbb{Z}}$, define a subsequence $\{k_n(\eta)\}_{k\in\mathbb{Z}}$ of $\mathbb{Z}$ given by
\[k_n(\eta)=\left\{\begin{array}{ll}
\min\left\{m\in\mathbb{Z}\::\:m>0,\ \eta_m\in\{0,i\}\right\}, & \mbox{if  }n=0,\\
\min\left\{m\in\mathbb{Z}\::\:m> k_{n-1}(\eta),\ \eta_m\in\{0,i\}\right\}, & \mbox{if  }n\geq1,\\
\max\left\{m\in\mathbb{Z}\::\:m< k_{n+1},\ \eta_m\in\{0,i\}\right\}, & \mbox{if  }n\leq-1,\\
\end{array}\right.\]
This $\{k_n(\eta)\}_{k\in\mathbb{Z}}$ is well-defined for $\eta$ almost everywhere.

For each $n\in\mathbb{Z}$, define transform $g_n:\{0,1,\cdots,\kappa\}^{\mathbb{Z}}\rightarrow\{0,i\}$ given by
\[g_n(\eta)=\eta_{k_n(\eta)}.\]
Then it holds that $\left(g_n(T_i\eta)\right)_{n\in\mathbb{Z}}=T_i\left(g_n(\eta)\right)_{n\in\mathbb{Z}}$, for each $\eta\in\{0,1,\cdots,\kappa\}^{\mathbb{Z}}$. Moreover, \cite{DS} shows that
\[T_i\left(g_n(\eta)\right)_{n\in\mathbb{Z}}\buildrel{d}\over{=}\left(g_n(\eta)\right)_{n\in\mathbb{Z}}.\]

Denote the filtration, 
\[\mathcal{F}:=\left\{f_n,n\in\mathbb{Z}\right\},\ \ \mathcal{G}:=\left\{g_n,n\in\mathbb{Z}\right\}.\]
It is obvious that $f_n(\eta)$ and $g_m(\eta)$ are independent for any $n,m\in\mathbb{Z}$, so is $\mathcal{F}$ and $\mathcal{G}$. Also $f_n(\eta)$ and $g_m(T_i\eta)$ are independent

Configuration $\eta$ is determined by $\left(f_n(\eta)\right)_{n\in\mathbb{Z}}$ and $\left(g_n(\eta)\right)_{n\in\mathbb{Z}}$, so there is a transform $\varphi$ such that 
\[\varphi\left(\left(f_n(\eta)\right)_{n\in\mathbb{Z}},\ \left(g_n(\eta)\right)_{n\in\mathbb{Z}}\right)=\eta,\ \ \forall \eta\in\{0,1,\cdots,\kappa\}^{\mathbb{Z}},\]
 which is measurable with respect to product measure $\mathcal{F}\times\mathcal{G}$.

Then it holds that
\begin{align*}
T_i\eta&=\varphi\left(\left(f_n(T_i\eta)\right)_{n\in\mathbb{Z}},\ \left(g_n(T_i\eta)\right)_{n\in\mathbb{Z}}\right)\\
&=\varphi\left(\left(f_n(\eta)\right)_{n\in\mathbb{Z}},\ T_i\left(\left(g_n(\eta)\right)_{n\in\mathbb{Z}}\right)\right)\\
&\buildrel{d}\over{=}\varphi\left(\left(f_n(\eta)\right)_{n\in\mathbb{Z}},\ \left(g_n(\eta)\right)_{n\in\mathbb{Z}}\right)\\
&=\eta.
\end{align*}
\end{proof}

{\cor As the same setting in Theorem \ref{iidinvariant}, it holds that
\[T\eta\buildrel{d}\over{=}\eta,\]
where $T=T_\kappa\circ\cdots\circ T_1$.
}

\subsection{Multicolor BBS on $\mathbb{R}$}\label{onR}

In this section, we consider a generalization of the multicolor BBS, whose dynamic is defined for continuous path in $\mathbb{R}^\kappa$. At first, we define Pitman transform for continuous path.

{\df \label{twoside pitman conti}
Let $\alpha\in{\mathbb{R}}^k,\ \alpha\ne0$. The two-sided Pitman transform $P_\alpha$ with respect to $\alpha$ is defined on the set \[\{\pi:\mathbb{R}\to {\mathbb{R}}^k,\ \pi(0)=0,\ \inf_{y\leq0}\alpha\cdot\pi(y)>-\infty\}\]by the formula, 
\[P_\alpha\pi(x)=\pi(x)-2\inf_{y\leq x}\frac{\alpha\cdot\pi(y)}{|\alpha|^2}\alpha+2\inf_{y\leq 0}\frac{\alpha\cdot\pi(y)}{|\alpha|^2}\alpha,\ \ \ x\in\mathbb{R} \]
Similarly to discrete case, for $k=1$, it holds that
\[P_\alpha\pi(x)=\pi(x)-2\inf_{y\leq x}\pi(y)+2\inf_{y\leq 0}\pi(y),\ \ \ x\in\mathbb{R}\]
for any $\alpha\in\mathbb{R},\ \alpha\ne0$, and it does not depend on $\alpha$. Then we define
\[P_1:=P_\alpha\ \ \ for\ \alpha \in\mathbb{R},\ \alpha\ne0.\]
Also we define the transform ${P_\alpha}^{-1}$ on the set \[\{\pi:\mathbb{R}\to {\mathbb{R}}^k,\ \pi(0)=0,\ \inf_{y\geq0}\alpha\cdot\pi(y)>-\infty\}\]by the formula, 
\[{P_\alpha}^{-1}\pi(x)=\pi(x)-2\inf_{y\geq x}\frac{\alpha\cdot\pi(y)}{|\alpha|^2}\alpha+2\inf_{y\geq 0}\frac{\alpha\cdot\pi(y)}{|\alpha|^2}\alpha,\ \ \ x\in\mathbb{R}, \]
and for $k=1$, 
\[{P_1}^{-1}\pi(x)=\pi(y)-2\inf_{y\geq x}\pi(y)+2\inf_{y\geq 0}\pi(y),\ \ \ x\in\mathbb{R}.\]
}

Unlike the discrete case, we can not describe the particle configuration $\eta$ directly, so we consider the dynamic for the path encoding S only. By analogy with the relevant discrete objects, define the path space

\[\mathcal{S}^0_c=\{S:\mathbb{R}\rightarrow\mathbb{R}^\kappa\::\:S_0=0,\ S \mbox{ is continuous} \}.\]

As the extension of \eqref{A2} in Remark \ref{Adef2} and \eqref{tau2} in Remark \ref{taudef2}, 
we define $A_i$ and $\tau$ as follows.
{\df Define $A_i:\mathcal{S}^0_c\rightarrow\mathbb{R}$ and $\tau_{(0,i)}:\mathcal{S}^0_c\rightarrow\mathcal{S}^0_c$ as follows\::
\[A_iS_x=-2\frac{(e_i-e_0)\cdot S_x}{|e_i-e_0|^2},\]
\[\tau_{(0,i)}S_x=S_x+A_iS_x(e_i-e_0)\]
for $x\in\mathbb{R}$.}

It is the case that 
the projection of $S_x$ along $e_i-e_0$ is $-\frac{1}{2}A_iS_x\left(e_i-e_0\right)$, and $S_x$ is decomposed into the sum as follows\::

\begin{equation}\label{projectionexpression}
S_x=-\frac12A_iS_x(e_i-e_0)+\left\{S_x+\frac12A_iS_x(e_i-e_0)\right\},
\end{equation}
and also it holds that
\begin{equation}\label{tauconti}
\tau_{(0,i)}S_x=\frac12A_iS_x(e_i-e_0)+\left\{S_x+\frac12A_iS_x(e_i-e_0)\right\}.
\end{equation}

\vspace{5pt}
Then we can define the dynamics of the generalized multicolor BBS, given by

\[T_i=\tau_{(0,i)}P_{e_i-e_0},\ \ \mbox{on}\ \{S\in\mathcal{S}^0_c\::\:\limsup_{x\rightarrow-\infty}A_iS_x<\infty\},\]

\[T^{-1}_i=P^{-1}_{e_i-e_0}\tau_{(0,i)},\ \ \mbox{on}\ \{S\in\mathcal{S}^0_c\::\:\liminf_{x\rightarrow\infty}A_iS_x>-\infty\}\]
for each $i\in\mathcal{C}$.

Moreover, the previous definitions of $A_i$ and $\tau_{(0,i)}$ yield the following alternative expression for $T_i$.
{\thm It holds that

\begin{equation}\label{Texpression}
T_iS_x=S_x+\left(A_iS_x-\sup_{y\leq x}A_iS_y+\sup_{y\leq 0}A_iS_y\right)(e_i-e_0),\ \ x\in\mathbb{R}
\end{equation}
for any $i\in\mathcal{C}$.
}

\begin{proof}
From \eqref{projectionexpression} and \eqref{tauconti}, it holds that
\begin{align*}
T_iS_x &=\tau_{(0,i)}P_{e_i-e_0}S_x\\
&=\tau_{(0,i)}P_{e_i-e_0}\left(-\frac12A_iS(e_i-e_0)+\left\{S+\frac12A_iS(e_i-e_0)\right\}\right)_x\\
&=\tau_{(0,i)}\left(P_{1}\left(-\frac12A_iS\right)_x(e_i-e_0)+\left\{S_x+\frac12A_iS_x(e_i-e_0)\right\}\right)\\
&=\tau_{(0,i)}\left(\frac12P_{1}\left(-A_iS\right)_x(e_i-e_0)+\left\{S_x+\frac12A_iS_x(e_i-e_0)\right\}\right)\\
&=-\frac12P_{1}\left(-A_iS\right)_x(e_i-e_0)+\left\{S_x+\frac12A_iS_x(e_i-e_0)\right\}\\
&=-\frac12\left\{-A_iS_x-2\inf_{y\leq x}(-A_iS_y)+2\inf_{y\leq 0}(-A_iS_y)\right\}(e_i-e_0)+\left\{S_x+\frac12A_iS_x(e_i-e_0)\right\}\\
&=S_x+\left(A_iS_x-\sup_{y\leq x}A_iS_y+\sup_{y\leq 0}A_iS_y\right)(e_i-e_0).
\end{align*}
\end{proof}

As in the discrete case, it is natural to seek to characterize the set

\[\mathcal{S}^{inv}_{\mathcal{C},c}:=\{S\in\mathcal{S}^0_c\::\:TS\in\bigcap_{i\in\mathcal{C}}\mathcal{S}^{rev}_{i,c}\mbox{\ for any composition }T\mbox{ of }T_i,\ i\in\mathcal{C}\},\]
where
\[\mathcal{S}^{rev}_{i,c}=\{S\in\mathcal{S}^0_c\::\:T_iS,T^{-1}_iS,T^{-1}_iTS,T_iT^{-1}_iS\mbox{ well-defined},\ T^{-1}_iT_iS=T_iT^{-1}_iS=S\}.\]
The following result is obtained by the similar argument in the discrete case. 
{\thm It holds that
\[\mathcal{S}^{inv}_{\mathcal{C},c}\supseteq\mathcal{S}^{good}_{\mathcal{C},c},\]
where 
\[\mathcal{S}^{good}_{\mathcal{C},c}:=\left\{S\in \mathcal{S}^0_c\::\:\forall i\in\mathcal{C}\ \exists F_i\in\mathcal{F}_c,\ \lim_{x\rightarrow\pm\infty}\frac{A_iS_x}{F_i(x)}=1\mbox{ and }\limsup_{x\rightarrow\pm\infty}\frac{F_j(x)}{F_i(x)}<\infty\ \forall i,j\in\mathcal{C}\right\},\]
and
\[\mathcal{F}_c=\{F:\mathbb{R}\rightarrow\mathbb{R}\::\:\mbox{increasing function, }\lim_{x\rightarrow\infty} F(x)=\infty,\ \lim_{x\rightarrow-\infty} F(x)=-\infty\}.\]

}

\subsection{Brownian motion with drift}\label{BMD}

Next, we consider a stochastic process whose path belongs to $\mathcal{S}^{good}_{\mathcal{C},c}$ almost surely. 
As an example, let $S=(S_x)_{x \in \mathbb{R}}$ be two-sided standard $\kappa-$dimensional standard Brownian motion with drift ${D}\in\mathbb{R}^\kappa$. Namely, for $x\geq 0$, we define $S_x=B^1_x+x{D}$, $S_{-x}=-\left(B^2_x+x{D}\right)$, where $B^1,B^2$ are independent standard Brownian motions in $\mathbb{R}^\kappa$. Since
\[A_iS_x=-2\frac{(e_i-e_0)\cdot B^1_x}{|e_i-e_0|^2}-2x\frac{(e_i-e_0)\cdot {D}}{|e_i-e_0|^2}\]
for $x\ge0$, the condition $\lim_{x\rightarrow\infty}\frac{A_iS_x}{F_i(x)}=1,\ a.s.\  \exists F_i\in\mathcal{F}_c$ is satisfied if and only if $(e_i-e_0)\cdot {D}<0$, and we can take $F_i(x)=-2x\frac{(e_i-e_0)\cdot {D}}{|e_i-e_0|^2},x\ge0$. Similarly, $\lim_{x\rightarrow-\infty}\frac{A_iS_x}{F_i(x)}=1,\ a.s.\  \exists F_i\in\mathcal{F}_c$ if and only if $(e_i-e_0)\cdot {D}<0$. Therefore, it holds that
\[S\in\mathcal{S}^{good}_{\mathcal{C},c},\ a.s.\ \Leftrightarrow\ (e_i-e_0)\cdot {D}<0,\ \ \forall i\in\mathcal{C}.\]
On the other hand, from Proposition \ref{vectorsproperty} $($\hspace{-1pt}ⅵ\hspace{-1pt}$)$
, there is an $\kappa+1$-tuple $c_0,\cdots,c_\kappa$ of real numbers for ${D}\in\mathbb{R}^\kappa$ such that
\[{D}=c_0e_0+\cdots+c_\kappa e_\kappa,\ \ c_0+\cdots+c_\kappa=0,\]
and, by Proposition \ref{vectorsproperty} $($\hspace{-1pt}ⅳ\hspace{-1pt}$)$, it holds that 
\begin{align*}
(e_i-e_0)\cdot{D}&=(e_i-e_0)\cdot\left(\frac12(c_i-c_0)(e_i-e_0)+\sum_{j\ne i}\left(c_j-\frac{c_i+c_0}{2}\right)e_j\right)\\
&=(c_i-c_0)\frac{|e_i-e_0|^2}{2}
\end{align*}
Thus we obtain the following set\::
\begin{align*}
\mathcal{D}:=&\{D\in\mathbb{R}^\kappa\::\:(e_i-e_0)\cdot {D}<0,\ \ \forall i\in\mathcal{C}\}\\
=&\{{D}\in\mathbb{R}^\kappa\::\:{D}=c_0e_0+\cdots+c_\kappa e_\kappa,\ c_0>c_i,\ \forall i\in\mathcal{C},\ c_0+\cdots+c_\kappa=0\},
\end{align*}
and it is the case that
\[S\in\mathcal{S}^{good}_{\mathcal{C},c},\ a.s.\ \Leftrightarrow\ D\in\mathcal{D}.\]

The main theorem in this subsection is the following which implies any Brownian motion with drift belonging to $\mathcal{S}^{good}_{\mathcal{C},c}$ is invariant under the actions of the generalized multicolor BBS. 
{\thm\label{BMinvariant} If $S$ is the two-sided $\kappa-$dimensional standard Brownian motion with drift ${D}\in\mathcal{D}$, then $T_iS \buildrel{d}\over{=}S$ for each $i\in\mathcal{C}$.}

{\cor As the same setting in Theorem \ref{BMinvariant}, it holds that
\[TS\buildrel{d}\over{=}S,\]
where $T=T_\kappa\circ\cdots\circ T_1$.
}

\vspace{10pt}

Before prove this main theorem, we show that Brownian motion with drift is obtained by a simple random walk scaling limit. 
From now on, fix $c_0,\cdots c_\kappa$ satisfying $c_0>c_i\ \forall i\in\mathcal{C},\ c_0+\cdots+c_\kappa=0$ and define
\[{D}=c_0e_0+\cdots+c_\kappa e_\kappa,\]
and
\[p^{(n)}_i=\frac1{\kappa+1}+\frac{c_i}{\sqrt {n\kappa}},\ \ i\in\mathcal{C}\cup\{0\}.\]
for large enough $n$ satisfying $0<p^{(n)}_i<1,\forall i\in\mathcal{C}$. 
Then we introduce vector valued random variables $\xi^{(n)}$  with distribution
\begin{equation}\label{xidist}
\mathbf{P}(\xi^{(n)}=e_i)=p^{(n)}_i,\ \ i\in\mathcal{C}\cup\{0\}.
\end{equation}
Moreover, let $\left\{\zeta^{(n)}_j\right\}_{j\in\mathbb{Z}}$ a sequence of independent identically distributed vector valued random variables and each $\zeta^{(n)}_j$ has the same distribution as $\xi^{(n)}$. Also we define the sequence of partial sums
\[S^{(n)}_{[x]}=\left\{\begin{array}{ll}
\zeta^{(n)}_1+\cdots+\zeta^{(n)}_{[x]}, & \mbox{if }[x]\ge1,\\
0, & \mbox{if }[x]=0,\\
-\left(\zeta^{(n)}_{-1}+\cdots+\zeta^{(n)}_{[x]}\right), & \mbox{if }[x]\le-1,
\end{array}\right.\]
and its linear interpolation
\begin{equation}\label{Y}
Y^{(n)}_x=S^{(n)}_{[x]}+\left(x-[x]\right)\zeta^{(n)}_{[x]+1},\ \ x\in\mathbb{R}.
\end{equation}

We introduce the notation $\mu^{p^{(n)}}$ to represent the probability measure on $\mathcal{S}^0_c$ induced by the stochastic process $(Y^{(n)}_x)_{x\in\mathbb{R}}$. As shown in theorem \ref{iidinvariant}, we have the invariance of $\mu^{p^{(n)}}$ under $T_i$ for any $i\in\mathcal{C}$.
As explained above, let $S=(S_x)_{x \in \mathbb{R}}$ be two-sided $\kappa-$dimensional Brownian motion with drift ${D}\in\mathbb{R}^\kappa$ and denote $\nu_D$ the probability measure on $\mathcal{S}^0_c$ induced by $S=(S_x)_{x \in \mathbb{R}}$.
Also we write $\mu_{a,b}$ to be the scaled measure such that
\[\mu_{a,b}\left(S \in A\right)= \mu\left( aS_{b \cdot} \in A\right),\]
for a probability measure $\mu$ on $\mathcal{S}^0_c$ and $a,b>0$.

The following theorem is known as the Invariance Principle of Donsker.
{\thm\label{donsker} $\nu_n:=\mu^{p^{(n)}}_{\frac{\sqrt\kappa}{\sqrt n},n}$ converges weakly to $\nu_{D}$.}

To prove this theorem, we prepare some lemmas.

{\lem \label{property5}
For any $u\in\mathbb{R}^\kappa$ satisfying $|u|=1$, it holds that
\[(e_0\cdot u)^2+\cdots(e_\kappa\cdot u)^2=\frac{\kappa+1}{\kappa}.\]
}
\begin{proof}
By Proposition \ref{vectorsproperty} $($\hspace{-1pt}ⅴ\hspace{-1pt}$)$, there are $a_0,\cdots,a_\kappa\in\mathbb{R}$ such that
\[u=a_0e_0+\cdots+a_\kappa e_\kappa\]
The condition $|u|=1$, \eqref{length} and \eqref{product} shows that
\[\sum^{\kappa}_{i=0}a^2_i-\frac2{\kappa}\sum_{i\ne j}a_ia_j=1\]
Then it holds that
\begin{align*}
\sum^{\kappa}_{i=0}(e_i\cdot u)^2&=\sum^{\kappa}_{i=0}\left(a_i-\frac1{\kappa}\sum_{j\ne i}a_j\right)^2\\
&=\sum^{\kappa}_{i=0}\left\{a^2_i-\frac{2a_i}{\kappa}\sum_{j\ne i}a_j+\frac1{\kappa^2}\left(\sum_{j\ne i}a_j\right)^2\right\}\\
&=\frac{\kappa+1}{\kappa}\sum^{\kappa}_{i=0}a^2_i-\frac{2(\kappa+1)}{\kappa^2}\sum_{i\ne j}a_ia_j\\
&=\frac{\kappa+1}{\kappa}.
\end{align*}
\end{proof}

{\lem \label{property6}
For any $u,\,v\in\mathbb{R}^\kappa$ satisfying $u\cdot v=0$, it holds that
\[(e_0\cdot u)(e_0\cdot v)+\cdots+(e_\kappa\cdot u)(e_\kappa\cdot v)=0.\]
}
\begin{proof}
By Proposition \ref{vectorsproperty} $($\hspace{-1pt}ⅴ\hspace{-1pt}$)$, there are $a_0,\cdots,a_\kappa,b_0,\cdots,b_\kappa\in\mathbb{R}$ such that
\[u=a_0e_0+\cdots+a_\kappa e_\kappa,\]
\[v=b_0e_0+\cdots+b_\kappa e_\kappa.\]
The condition $u\cdot v=0$, \eqref{length} and \eqref{product} show that
\[\sum^{\kappa}_{i=0}a_ib_i-\frac1{\kappa}\sum_{i\ne j}a_ib_j=0.\]
Then it holds that
\begin{align*}
\sum^{\kappa}_{i=0}(e_i\cdot u)(e_i\cdot v)&=\sum^{\kappa}_{i=0}\left(a_i-\frac1{\kappa}\sum_{j\ne i}a_j\right)\left(b_i-\frac1{\kappa}\sum_{j\ne i}b_j\right)\\
&=\frac{\kappa+1}{\kappa}\sum^{\kappa}_{i=0}a_ib_i-\frac{\kappa+1}{\kappa^2}\sum_{i\ne j}a_ib_j=0.
\end{align*}
\end{proof}

{\lem \label{property8}
$(a)$\ For each $i\in\mathcal{C}\cup\{0\}$, we denote the components of $e_i$ as follows\::

\[e_0=\left(\begin{array}{c}e_{0,1}\\e_{0,2}\\e_{0,3}\\ \vdots\\ e_{0,\kappa-1}\\e_{0,\kappa}\end{array}\right),\ \ 
e_1=\left(\begin{array}{c}e_{1,1}\\e_{1,2}\\e_{1,3}\\ \vdots\\ e_{1,\kappa-1}\\e_{1,\kappa}\end{array}\right),\ \cdots,\ 
e_\kappa=\left(\begin{array}{c}e_{\kappa,1}\\e_{\kappa,2}\\e_{\kappa,3}\\ \vdots\\ e_{\kappa,\kappa-1}\\e_{\kappa,\kappa}\end{array}\right).\]
For any $s,t\in\mathcal{C},s\ne t$ it holds that
\[e^2_{0,s}+\cdots+e^2_{\kappa,s}=\frac{\kappa+1}{\kappa},\]
\[e_{0,s}e_{0,t}+\cdots+e_{\kappa,s}e_{\kappa,t}=0.\]
$(b)$\ For each $n$, we denote the components of $\xi^{(n)}$ by $\xi^{(n)}=(\xi^{(n)}_{1},\cdots,\xi^{(n)}_{\kappa})$.
For any $s,t\in\mathcal{C},s\ne t$ it holds that
\[\lim_{n\rightarrow\infty}\mathbf{E}(\xi^{(n)}_s)=0,\ \lim_{n\rightarrow\infty}\mathbf{V}(\xi^{(n)}_s)=\frac1{\kappa},\ 
\lim_{n\rightarrow\infty}\mathbf{E}(\xi^{(n)}_s\xi^{(n)}_s)=0.\]
}

\begin{proof}
In Lemma \ref{property5} and \ref{property6}, let $u=(\delta_{s\,1},\cdots,\delta_{s\,\kappa})$ and $v=(\delta_{t\,1},\cdots,\delta_{t\,\kappa})$ for $s,t\in\mathcal{C},s\ne t$, where $\delta$ is the Kronecker delta. Then the two equations in (a) follow directly.

Assume that $\xi$ is vector valued random variable with distribution
\begin{equation}\label{uniformprob}
\mathbf{P}(\xi=e_i)=\frac{1}{\kappa+1},\ \ i\in\mathcal{C}\cup\{0\},
\end{equation}
and denote its components by $\xi=(\xi_1,\cdots,\xi_\kappa)$. Then, by Proposition \ref{vectorsproperty} , 
\[\mathbf{E}(\xi_s)=0,\ \ s\in\mathcal{C},\]
where $\mathbf{E}$ is the expectation with respect to $\mathbf{P}$. Also above equations in (a) show that
\[\mathbf{V}(\xi_s)=\frac1{\kappa},\ \ s\in\mathcal{C}\]
where $\mathbf{V}$ is the variance with respect to $\mathbf{P}$, and
\[\mathbf{E}(\xi_s\xi_t)=0,\ \ s,t\in\mathcal{C},s\ne t.\]
The distribution \eqref{xidist} and \eqref{uniformprob} imply that $\xi^{(n)}$ converges to $\xi$ almost surely as $n\rightarrow\infty$, and convergence theorem shows the claim (b).
\end{proof}

{\rem \label{Dcomponent} Denote the components of $D=c_0e_0+\cdots+c_\kappa e_\kappa$ by $D=(D_1,\cdots,D_\kappa)$. Then it holds that
\[E\left(\xi^{(n)}_j\right)=\frac{1}{\sqrt {n\kappa}}D_j,\ \ j\in\mathcal{C}.\]
It follows directly by \eqref{xidist}.
}

\vspace{10pt}
To prove Theorem \ref{donsker}, it is enough to show following two claims.\\
(1) The finite-dimensional distribution of $\nu_n$ converges weakly to that of $\nu_{D}$.\\
(2) $\left\{\nu_n\right\}_n$ is tight.

\vspace{5pt}
We prove (1) as Proposition \ref{finitedist} and show what is equivalent to (2) as Proposition \ref{tight}. In the proof of Proposition \ref{finitedist} and \ref{tight}, we write $|\cdot|$ as the Euclidean norm.

{\prop\label{finitedist}
Define the stochastic process
\[X^{(n)}_x=\frac{\sqrt\kappa}{\sqrt n}Y^{(n)}_{nx},\]
where $Y^{(n)}$ is given by \eqref{Y}. Then, for any $0\leq x_1<\cdots<x_d<\infty$, 
\[\left(X^{(n)}_{x_1},\cdots,X^{(n)}_{x_d}\right)\buildrel{d}\over{\rightarrow}\left(B_{x_1}+x_1{D},\cdots,B_{x_d}+x_d {D}\right)\ \ \ \ \mbox{as}\ n\rightarrow\infty\]
where $\{B_x\}_{x\ge0}$ is a $\kappa$-dimensional Brownian motion. Also the same is true for $x\le0$.}

\begin{proof}

We prove the case $d=2$, that is 
\[\left(X^{(n)}_{s},X^{(n)}_{t}\right)\buildrel{d}\over{\rightarrow}\left(B_{s}+s{D},B_{t}+t{D}\right)\ \ \mbox{for}\ 0<s<t,\]
and the other case proved samely. Since 
\[\left|X^{(n)}_x-\frac{\sqrt\kappa}{\sqrt n}S^{(n)}_{[nx]}\right| \leq \frac{\sqrt\kappa}{\sqrt n}\left|\zeta^{(n)}_{[nx]+1}\right|=\frac{\sqrt\kappa}{\sqrt n},\]
we have by the Chebyshev inequality,
\[\mathbf{P}\left(\left|X^{(n)}_x-\frac{\sqrt\kappa}{\sqrt n}S^{(n)}_{[nx]}\right|>\varepsilon\right) \leq \frac{\kappa}{\varepsilon^2n}\rightarrow0\]
as $n\rightarrow\infty$. Then it is clear that 
\[\left|\left(X^{(n)}_{s},X^{(n)}_{t}\right)-\frac{\sqrt\kappa}{\sqrt n}\left(S^{(n)}_{[ns]},S^{(n)}_{[nt]}\right)\right|\rightarrow0\ \ \ \mbox{in probability.}\]
Therefore, it is enough to show that 
\[\frac{\sqrt\kappa}{\sqrt n}\left(S^{(n)}_{[ns]},S^{(n)}_{[nt]}\right)\buildrel{d}\over{\rightarrow}\left(B_{s}+s{D},B_{t}+t{D}\right), \]
and it is equivalent to
\[\frac{\sqrt\kappa}{\sqrt n}\left(\sum^{[ns]}_{m=1}\zeta^{(n)}_m,\sum^{[nt]}_{m=[ns]+1}\zeta^{(n)}_m\right)\buildrel{d}\over{\rightarrow}\left(B_{s}+s{D},B_{t}-B_{s}+(t-s) {D}\right).\]
The independence of the random variables $\{\zeta_m\}^{\infty}_{m=1}$ implies
\begin{align*}
&\mathbf{E}\left(\exp\left\{\frac{\sqrt\kappa}{\sqrt n}i\left(\sum^{[ns]}_{m=1}\zeta^{(n)}_m\cdot u+\sum^{[nt]}_{m=[ns]+1}\zeta^{(n)}_m\cdot v\right)\right\}\right)\\
={}&\mathbf{E}\left(\exp\left\{\frac{\sqrt\kappa}{\sqrt n}i\sum^{[ns]}_{m=1}\zeta^{(n)}_m\cdot u\right\}\right)\mathbf{E}\left(\exp\left\{\frac{\sqrt\kappa}{\sqrt n}i\sum^{[nt]}_{m=[ns]+1}\zeta^{(n)}_m\cdot v\right\}\right),
\end{align*}
for any $u=(u_1,\cdots,u_\kappa),v=(v_1,\cdots,v_\kappa)\in\mathbb{R}^\kappa$, and also it holds that
\[\mathbf{E}\left(\exp\left\{\frac{\sqrt\kappa}{\sqrt n}i\sum^{[ns]}_{m=1}\zeta^{(n)}_m\cdot u\right\}\right)=
\left\{\varphi_n\left(\frac{\sqrt\kappa}{\sqrt n}u\right)\right\}^{[ns]},\]
where $\varphi_n(\theta)$ is the characteristic function of $\xi^{(n)}$ given by
\[\varphi_n(\theta)=\mathbf{E}\left(\exp\left\{i\xi^{(n)}\cdot\theta\right\}\right)\]
for $\theta=(\theta_j)_{1\le j\le\kappa}\in\mathbb{R}^\kappa$. The function $\varphi_n$ satisfies 
\begin{align*}
\frac{\partial \varphi_n}{\partial \theta_j}(\theta)&=\mathbf{E}\left(i\xi^{(n)}_j\exp\left\{i\xi^{(n)}\cdot\theta\right\}\right)\ \ j\in\mathcal{C},\\
\frac{\partial \varphi_n}{\partial^2 \theta_j}(\theta)&=\mathbf{E}\left(-{\xi^{(n)}_j}^2\exp\left\{i\xi^{(n)}\cdot\theta\right\}\right)\ \ j\in\mathcal{C},\\
\frac{\partial \varphi_n}{\partial \theta_j\partial \theta_k}(\theta)&=\mathbf{E}\left(-\xi^{(n)}_j\xi^{(n)}_k\exp\left\{i\xi^{(n)}\cdot\theta\right\}\right)\ \ j,k\in\mathcal{C}, j\ne k.
\end{align*}
Remark \ref{Dcomponent} implies
\[\frac{\partial \varphi_n}{\partial \theta_j}(0)=\frac{iD_j}{\sqrt {n\kappa}},\]
and Lemma \ref{property8} shows, if $n\rightarrow\infty$ and $\theta\rightarrow0$, that is, $\theta_j\rightarrow0$ for any $j\in\mathcal{C}$,
\[\frac{\partial \varphi_n}{\partial^2 \theta_j}(\theta)\rightarrow-\frac1{\kappa}\ \ \ 
\frac{\partial \varphi_n}{\partial \theta_j\partial \theta_k}(\theta)\rightarrow0.\]
By Taylor's theorem, there is a vector $u'=(u'_1,\cdots,u'_\kappa)\in\mathbb{R}^\kappa$ for fixed $u\in\mathbb{R}^\kappa$, such that $0\le u'_j\le \frac{\sqrt\kappa}{\sqrt n}u_j$ for any $j\in\mathcal{C}$ and 
\begin{align*}
&\varphi_n\left(\frac{\sqrt\kappa}{\sqrt n}u\right)\\
={}&\varphi_n(0)+\sum_{j\in\mathcal{C}}\frac{\sqrt\kappa}{\sqrt n}u_j\frac{\partial \varphi_n}{\partial \theta_j}(0)+\sum_{j\in\mathcal{C}}\frac12\left(\frac{\sqrt\kappa}{\sqrt n}u_j\right)^2\frac{\partial \varphi_n}{\partial^2 \theta_j}(u')+\sum_{j\ne k}\frac12\left(\frac{\kappa}{n}u_ju_k\right)\frac{\partial \varphi_n}{\partial \theta_j\partial \theta_k}(u')\\
={}&1+\sum_{j\in\mathcal{C}}\frac{\sqrt\kappa}{\sqrt n}u_j\frac{iD_j}{\sqrt {n\kappa}}+\sum_{j\in\mathcal{C}}\frac12\left(\frac{\sqrt\kappa}{\sqrt n}u_j\right)^2\frac{\partial \varphi_n}{\partial^2 \theta_j}(u')+\sum_{j\ne k}\frac12\left(\frac{\kappa}{n}u_ju_k\right)\frac{\partial \varphi_n}{\partial \theta_j\partial \theta_k}(u').
\end{align*}
Since $\log(1+x)=x+o(x)$ as $x\rightarrow0$, it holds that

\begin{align*}
&\log\left\{\varphi_n\left(\frac{\sqrt\kappa}{\sqrt n}u\right)\right\}^{[ns]}\\
={}&[ns]\log\varphi_n\left(\frac{\sqrt\kappa}{\sqrt n}u\right)\\
={}&i\frac{[ns]}{n}\sum_{j\in\mathcal{C}}D_ju_j+\frac{[ns]}{2n}\sum_{j\in\mathcal{C}}u^2_j\kappa\frac{\partial \varphi_n}{\partial^2 \theta_j}(u')+\frac{[ns]\kappa}{n}\sum_{j\ne k}\frac12\left(u_ju_k\right)\frac{\partial \varphi_n}{\partial \theta_j\partial \theta_k}(u')\\
\rightarrow{}&is {D}\cdot u-\sum_{j\in\mathcal{C}}\frac{su^2_j}{2}
\end{align*}
as $n\rightarrow\infty$. Thus,

\[\lim_{n\rightarrow\infty}\mathbf{E}\left(\exp\left\{\frac{\sqrt\kappa}{\sqrt n}i\sum^{[ns]}_{m=1}\zeta^{(n)}_m\cdot u\right\}\right)=\exp\left\{is {D}\cdot u-\sum_{j\in\mathcal{C}}\frac{su^2_j}{2}\right\}.\]
Similarly, 
\[\lim_{n\rightarrow\infty}\mathbf{E}\left(\exp\left\{\frac{\sqrt\kappa}{\sqrt n}i\sum^{[nt]}_{m=[ns]+1}\zeta^{(n)}_m\cdot v\right\}\right)=\exp\left\{i(t-s) {D}\cdot v-\sum_{j\in\mathcal{C}}\frac{(t-s)v^2_j}{2}\right\},\]
and the proof is complete.
\end{proof}

The tightness of $\{\nu_n\}_n$ is known to be equivalent to the following proposition \cite[Theorem.2.4.10, 2.4.15]{KS}.
{\prop\label{tight} With the same setting in Proposition \ref{finitedist}, it holds that

\begin{equation}\label{tight1}
\lim_{\lambda\uparrow\infty}\sup_{n\ge1}\mathbf{P}\left(\left|X^{(n)}_0\right|>\lambda\right)=0,
\end{equation}
and, for any $T>0,\varepsilon>0$, 
\begin{equation}\label{tight2}
\lim_{\delta\downarrow0}\max_{n\ge1}\mathbf{P}\left(\sup_{|t-s|\le \delta\ 0\le s,t\le T}\left|X^{(n)}_t-X^{(n)}_s\right|>\varepsilon\right)=0.
\end{equation}

}
\begin{proof}
Since $X^{(n)}_0=0$ for every $n$, \eqref{tight1} is obvious. We may replace $\sup_{n\ge1}$ in \eqref{tight2} by $\limsup_{n\rightarrow\infty}$ because for a finite number of integers $n$ we can make the probability appearing in \eqref{tight2} as small as we choose by reducing $\delta$. Let $X^{(n)}_t=\left(X^{(n)}_{t,1},\cdots,X^{(n)}_{t,\kappa}\right)$ for $t\geq0$, it holds that
\begin{align*}
\max_{|t-s|\le \delta\ 0\le s,t\le T}\left|X^{(n)}_t-X^{(n)}_s\right|&=\max_{|t-s|\le \delta\ 0\le s,t\le T}\sqrt{\sum^{\kappa}_{j=1}\left|X^{(n)}_{t,j}-X^{(n)}_{t,j}\right|^2}\\
&\leq \kappa \sum^{\kappa}_{j=1}\max_{|t-s|\le \delta\ 0\le s,t\le T}\left|X^{(n)}_{t,j}-X^{(n)}_{t,j}\right|.
\end{align*}
Thus, 
\begin{align*}
&\lim_{\delta\downarrow0}\limsup_{n\rightarrow\infty}\mathbf{P}\left(\max_{|t-s|\le \delta\ 0\le s,t\le T}\left|X^{(n)}_t-X^{(n)}_s\right|>\varepsilon\right)\\
\le{}&\lim_{\delta\downarrow0}\limsup_{n\rightarrow\infty}\mathbf{P}\left(\bigcup_{j\in\mathcal{C}}\max_{|t-s|\le \delta\ 0\le s,t\le T}\left|X^{(n)}_{t,j}-X^{(n)}_{s,j}\right|>\frac{\varepsilon}{\kappa}\right)\\
\le{}&\sum^{\kappa}_{j=1}\lim_{\delta\downarrow0}\limsup_{n\rightarrow\infty}\mathbf{P}\left(\max_{|t-s|\le \delta\ 0\le s,t\le T}\left|X^{(n)}_{t,j}-X^{(n)}_{s,j}\right|>\frac{\varepsilon}{\kappa}\right).
\end{align*}
By the definition of $X^{(n)}$, $Y^{(n)}$ and $S^{(n)}$, it holds that 
\[\mathbf{P}\left(\max_{|t-s|\le \delta\ 0\le s,t\le T}\left|X^{(n)}_{t,j}-X^{(n)}_{s,j}\right|>\frac{\varepsilon}{\kappa}\right)
=\mathbf{P}\left(\max_{|t-s|\le n\delta\ 0\le s,t\le nT}\left|Y^{(n)}_{t,j}-Y^{(n)}_{s,j}\right|>\frac{\varepsilon\sqrt n}{\kappa\sqrt\kappa}\right),\]
and
\begin{align*}
\max_{|t-s|\le n\delta\ 0\le s,t\le nT}\left|Y^{(n)}_{t,j}-Y^{(n)}_{s,j}\right|&\le\max_{|t-s|\le [n\delta]+1\ 0\le s,t\le [nT]+1}\left|Y^{(n)}_{t,j}-Y^{(n)}_{s,j}\right|\\
&\le\max_{1\le m\le[n\delta]+1\ 0\le k\le [nT]+1}\left|S^{(n)}_{m+k,j}-S^{(n)}_{k,j}\right|,
\end{align*}
where $Y^{(n)}_x=(Y^{(n)}_{x,1},\cdots,Y^{(n)}_{x,\kappa})$ and $S^{(n)}_m=(S^{(n)}_{m,1},\cdots,S^{(n)}_{m,\kappa})$. Therefore it is enough to show that
\[\lim_{\delta\downarrow0}\limsup_{n\rightarrow\infty}\mathbf{P}\left(\max_{1\le m\le[n\delta]+1\ 0\le k\le [nT]+1}\left|S^{(n)}_{m+k,j}-S^{(n)}_{k,j}\right|>\frac{\varepsilon\sqrt n}{\kappa\sqrt\kappa}\right)=0\]
for each $j\in\mathcal{C}$.

Recall the definition of $S^{(n)}$, the $j$-th component of it is as follows

\[S^{(n)}_{0,j}=0,\ S^{(n)}_{m,j}=\zeta^{(n)}_{1,j}+\cdots+\zeta^{(n)}_{m,j},\ \ m\ge1,\]
where $\{\zeta^{(n)}_{\ell,j}\}_{\ell\ge1}$ are independent and 
\begin{align*}
\mathbf{P}\left(\zeta^{(n)}_{\ell,j}=e_{i,j}\right)&=\frac{1}{\kappa+1}+\frac{c_i}{\sqrt {n\kappa}},\ \ i\in\mathcal{C},\\
\mathbf{E}\left(\zeta^{(n)}_{\ell,j}\right)&=\frac{D_j}{\sqrt {n\kappa}},
\end{align*}
from Remark \ref{Dcomponent}. Now we define a new stochastic process, 
\[R^{(n)}_{0,j}=0,\ R^{(n)}_{m,j}=\left(\zeta^{(n)}_{1,j}-\frac{D_j}{\sqrt {n\kappa}}\right)+\cdots+\left(\zeta^{(n)}_{m,j}-\frac{D_j}{\sqrt {n\kappa}}\right),\ \ m\ge1.\]
Since
\begin{align*}
&\left|S^{(n)}_{m+k,j}-S^{(n)}_{k,j}\right|\\
={}&\left|\left(R^{(n)}_{m+k,j}+(m+k)\frac{D_j}{\sqrt {n\kappa}}\right)-\left(R^{(n)}_{k,j}+k\frac{D_j}{\sqrt {n\kappa}}\right)\right|\\
\le{}&\left|R^{(n)}_{m+k,j}-R^{(n)}_{k,j}\right|+\left|m\frac{D_j}{\sqrt {n\kappa}}\right|,
\end{align*}
it is enough to show 
\[\lim_{\delta\downarrow0}\limsup_{n\rightarrow\infty}\mathbf{P}\left(\max_{1\le m\le[n\delta]+1\ 0\le k\le [nT]+1}\left|R^{(n)}_{m+k,j}-R^{(n)}_{k,j}\right|>\frac{\varepsilon\sqrt n}{2\kappa\sqrt\kappa}\right)=0,\]
and 
\[\lim_{\delta\downarrow0}\limsup_{n\rightarrow\infty}\mathbf{P}\left(\max_{1\le m\le[n\delta]+1}\left|\frac{m}{\sqrt n}D_j\right|>\frac{\varepsilon\sqrt {n\kappa}}{2\kappa\sqrt\kappa}\right)=0.\]
The first one is shown in \cite[Lemma.2.4.19]{KS}. Also it holds that

\begin{align*}
&\limsup_{n\rightarrow\infty}\mathbf{P}\left(\max_{1\le m\le[n\delta]+1}\left|\frac{m}{\sqrt {n\kappa}}D_j\right|>\frac{\varepsilon\sqrt n}{2\kappa\sqrt\kappa}\right)\\
={}&\limsup_{n\rightarrow\infty}\mathbf{P}\left(\left|\frac{[n\delta]+1}{\sqrt {n\kappa}}D_j\right|>\frac{\varepsilon\sqrt n}{2\kappa\sqrt\kappa}\right)=0,\ \ \mbox{if}\ \delta<\frac{\varepsilon}{2\kappa|D_j|}
\end{align*}
and this shows the second one. 
\end{proof}

Proposition \ref{finitedist} and \ref{tight} show Theorem \ref{donsker}.

\vspace{10pt}
Next, to prove Theorem \ref{BMinvariant}, we show following lemmas.

{\lem \label{ab}Let $a, b>0,i\in\mathcal{C}$. If $\mu$ is invariant under $T_i$, then $\mu_{a,b}$ is also invariant under $T_i$.}

\begin{proof}
Let $S^{a,b}_x=aS_{bx}$ for $a,b >0$ and $x \in \mathbb{R}$. By using the expression \eqref{Texpression}, it holds that 
\begin{align*}
T_iS^{a,b}_x&=aS_{bx}+\left(A_iaS_{bx}-\sup_{y\leq x}A_iaS_{by}+\sup_{y\leq 0}A_iaS_{by}\right)(e_i-e_0)\\
&=(T_iS)^{a,b}_x,
\end{align*}
and the claim follows.
\end{proof}

{\lem\label{generaltheory}
Suppose $\{\mu_n\}$ is a sequence of probability measures on $\mathcal{S}^0_c$, each of which is invariant under $T_i$, and $\mu_n$ converges weakly to $\mu$. Moreover, suppose that $\mu_n$ satisfies for any $z \in \mathbb{R}$,
\[\lim_{x \to -\infty} \limsup_{n \to \infty} \mu_n\left(\sup_{y\leq x}A_iS_{y} > A_iS_z\right)=0\]
and $\mu$ satisfies for any $z \in \mathbb{R}$,
\[\lim_{x \to -\infty} \mu\left(\sup_{y\leq x}A_iS_{y} > A_iS_z\right)=0.\]
It then holds that $\mu$ is also invariant under $T_i$.
}

\begin{proof}
It is enough to show that for any $L>0$ and continuous bounded function $f: C([-L,L],\mathbb{R}^\kappa) \to \mathbb{R}$,
\[\mu\left(f \left(S|_{[-L,L]}\right)\right) = \mu\left(f \left(T_iS|_{[-L,L]}\right)\right).\]
Let
\[
M^{L'}_x:=\left\{\begin{array}{ll}
       A_iS_{-L'}, & \mbox{if }x < -L', \\
    \sup_{-L' \le y \le x} A_iS_y, & \mbox{if }-L' \le x \le L',\\
       \sup_{-L' \le y \le L'} A_iS_y, & \mbox{otherwise.}
     \end{array}\right.\]
Also, define 

\begin{equation}\label{413}
(T^{L'}_iS)_x:=S_x+\left(A_iS_x-M^{L'}_x+M^{L'}_0\right)(e_i-e_0),\ \ x\in\mathbb{R}.
\end{equation}
Then, $T^{L'}_i :\mathcal{S}^0_c \to \mathcal{S}^0_c $ is continuous, and so
\begin{equation}\label{413a}
\lim_{n \to \infty} \mu_n \left(f \left((T^{L'}_iS)|_{[-L,L]}\right)\right) = \mu\left(f \left((T^{L'}_iS)|_{[-L,L]}\right)\right),
\end{equation}
for any $L,L'$. 

It is easy to verify that $(T^{L'}_iS)|_{[-L,L]}=(T_iS)|_{[-L,L]}$ if $L < L'$ and $\sup_{y\leq -L'}A_iS_{y} \leq A_iS_{-L}$, by comparing \eqref{Texpression} and \eqref{413}. Therefore, for any $L' >L$,
\[\left|\mu_n\left(f \left((T^{L'}_iS)|_{[-L,L]}\right)\right) - \mu_n\left(f \left((T_iS)|_{[-L,L]}\right)\right)\right| \le 2 \|f\|_{\infty} \mu_n \left(\sup_{y\leq -L'}A_iS_{y} > A_iS_{-L}\right).\]
Hence, by assumption, we have that
\[\lim_{L' \to \infty} \limsup_{n \to \infty}\left|\mu_n\left(f \left((T^{L'}S)|_{[-L,L]}\right)\right) - \mu_n\left(f \left((TS)|_{[-L,L]}\right)\right)\right|=0,\]
which implies, with \eqref{413a}, 
\begin{equation}\label{413aa}
\limsup_{n \to \infty}\mu_n\left(f \left((T_iS)|_{[-L,L]}\right)\right)=\lim_{L' \to \infty}  \mu(f (T^{L'}_iS|_{[-L,L]})).
\end{equation}
Similarly it holds that
\[\left|\mu\left(f \left((T^{L'}_iS)|_{[-L,L]}\right)\right) - \mu\left(f \left((T_iS)|_{[-L,L]}\right)\right)\right| \le 2 \|f\|_{\infty} \mu \left(\sup_{y\leq -L'}A_iS_{y} > A_iS_{-L}\right),\]
and the assumption $\lim_{x \to -\infty} \mu\left(\sup_{y\leq x}A_iS_{y} > A_iS_z\right)=0$ for any $x$ implies that 
\[\lim_{L' \to \infty} \mu(f (T^{L'}_iS|_{[-L,L]}))=\mu\left(f \left((T_iS)|_{[-L,L]}\right)\right).\]
This is the right-hand side of \eqref{413aa}. Also the left-hand side of \eqref{413aa} is equal to
\[\limsup_{n \to \infty}\mu_n\left(f \left(S|_{[-L,L]}\right)\right)=\mu\left(f \left(S|_{[-L,L]}\right)\right),\]then the claim is proved.
\end{proof}

Finally, We check the assumptions of the previous result for $\nu_n=\mu^{p^{(n)}}_{\frac{\sqrt\kappa}{\sqrt n},n}$ and $\nu_D$.

{\lem\label{checkassumption} For any $z \in \mathbb{R}$,
\[\lim_{x \to -\infty} \limsup_{n \to \infty} \nu_n\left(\sup_{y\leq x}A_iS_{y} > A_iS_z\right)=0\]
and
\[\lim_{x \to -\infty} \nu_D\left(\sup_{y\leq x}A_iS_{y} > A_iS_z\right)=0.\]
}

\begin{proof} 
For $S_x=B_x+(c_0e_0+\cdots c_\kappa e_\kappa)x$ it holds that
\begin{align*}
A_iS_x&=-2\frac{(e_i-e_0)\cdot S_x}{|e_i-e_0|^2}\\
&=-2\frac{(e_i-e_0)\cdot B_x}{|e_i-e_0|^2}+(c_0-c_i)x\\
&\rightarrow-\infty
\end{align*}
almost surely as $x\rightarrow-\infty$. Therefore, the second claim of the lemma is obvious. 
To estimate the probability $\nu_n\left(\sup_{y\leq x}A_iS_{y} > A_iS_z\right)$, first note that, for any $x<z$,
\begin{align*}
\nu_n \left(\sup_{y\leq x}A_iS_{y} > A_iS_z\right) & =\mu^{p^{(n)}}\left(\sup_{y\leq nx}A_iS_{y} > A_iS_{nz}\right)\\
& \leq \mu^{p^{(n)}} \left(\sup_{y\leq [nx]+1}A_iS_{y}> \min\left\{A_iS_{[nz]},A_iS_{[nz]+1} \right\}\right)\\
& = \mu^{p^{(n)}} \left(\sup_{y\leq [nx]+1-[nz]}A_iS_{y}> \min\left\{A_iS_{0},A_iS_{1} \right\}\right)\\
& \le \mu^{p^{(n)}} \left(\sup_{y\leq [nx]+1-[nz]}A_iS_{y}\ge 0 \right),
\end{align*}
where $[w]$ is the maximum integer not greater than $w$. Thus we only need to show that
\[\lim_{x \to -\infty} \limsup_{n \to \infty} \mu^{p^{(n)}} \left(\sup_{y\leq [nx]}A_iS_{y} \ge 0 \right)=0.\]
For any $\ell \ge 1$, we have
\begin{align*}
\mu^{p^{(n)}} \left(\sup_{y\leq -\ell}A_iS_{y} \ge 0 \right) &\leq  \mu^{p^{(n)}} \left(A_iS_{-\ell} \ge 0\right) + \sum_{k  \le -1}\mu^{p^{(n)}} \left(A_iS_{-\ell}=k\right)\left(\frac{p_i}{ p_0}\right)^{-k}.
\end{align*}
Now, since $A_iS_{-\ell}\buildrel{d}\over{=}-A_iS_\ell=-\sum_{k=1}^{\ell}\left(\mathbf{1}_{\{\eta_k=0\}}-\mathbf{1}_{\{\eta_k=i\}}\right)$, we have
\begin{align*}
\mu^{p^{(n)}} \left(A_iS_{-\ell} \ge 0\right) &= \mu^{p^{(n)}}\left(\sum_{k=1}^{\ell}\left(\mathbf{1}_{\{\eta_k=0\}}-\mathbf{1}_{\{\eta_k=i\}}\right) \le 0\right)\\
& \le \mu^{p^{(n)}}\left(\frac{1}{\ell}\left|\sum_{k=1}^{\ell}\left(\mathbf{1}_{\{\eta_k=0\}}-\mathbf{1}_{\{\eta_k=i\}}\right) - \frac{\ell (c_0-c_i)}{\sqrt {n\kappa}}\right| \ge \frac{c_0-c_i}{\sqrt {n\kappa}}\right)\\
& \le \mu^{p^{(n)}}\left(\frac{1}{\ell}\left|\sum_{k=1}^{\ell}\left\{\mathbf{1}_{\{\eta_k=0\}}-\mathbf{1}_{\{\eta_k=i\}}- \frac{c_0-c_i}{\sqrt {n\kappa}}\right\} \right| \ge \frac{c_0-c_i}{\sqrt {n\kappa}}\right)\\
& \le \frac{n\kappa}{\ell (c_0-c_i)^2} E^{p^{(n)}}\left(\left\{\mathbf{1}_{\{\eta_k=0\}}-\mathbf{1}_{\{\eta_k=i\}}-\frac{c_0-c_i}{\sqrt {n\kappa}}\right\}^2\right)\\
&= \frac{n\kappa}{\ell (c_0-c_i)^2} V^{p^{(n)}}\left(\mathbf{1}_{\{\eta_k=0\}}-\mathbf{1}_{\{\eta_k=i\}}\right)\\
&=\frac{n\kappa}{\ell (c_0-c_i)^2}\left(\frac{2}{\kappa+1}+\frac{c_0+c_i}{\sqrt {n\kappa}}\right),
\end{align*}
where $E^{p^{(n)}}$ is the expectation and $V^{p^{(n)}}$ is the variance with respect to $\mu^{p^{(n)}}$.  Moreover, 
\begin{align*}
 &\sum_{k  \le -1}\mu^{p^{(n)}} \left(A_iS_{-\ell}=k\right)\left(\frac{p_i}{ p_0}\right)^{-k}\\
 ={}&\sum_{-\ell\leq k  \le -1}\mu^{p^{(n)}} \left(A_iS_{-\ell}=k\right)\left(\frac{p_i}{ p_0}\right)^{-k}\\
  ={}&\sum_{-\ell\leq k  \le -1}\sum_{0\le j\le \ell+k}\binom{\ell}{j}(1-p_0-p_i)^{j}\binom{\ell-j}{\frac{\ell-j-k}{2}}p^{\frac{\ell-j-k}{2}}_0p^{\frac{\ell-j+k}{2}}_i\left(\frac{p_i}{ p_0}\right)^{-k}\\
 ={}&\sum_{-\ell\leq k  \le -1}\sum_{0\le j\le \ell+k}\binom{\ell}{j}(1-p_0-p_i)^{j}\binom{\ell-j}{\frac{\ell-j-k}{2}}p^{\frac{\ell-j+k}{2}}_0p^{\frac{\ell-j-k}{2}}_i\\
  ={}&\sum_{-\ell\leq k  \le -1}\mu^{p^{(n)}} \left(A_iS_{-\ell}=-k\right)\\
  ={}&\mu^{p^{(n)}} (A_iS_{-\ell} \ge 1)\\
  \le{}&\mu^{p^{(n)}} (A_iS_{-\ell} \ge 0),\\
\end{align*}
where $\binom{\ell}{q} \equiv 0$ for $q \notin \mathbb{N}$. Therefore, we have
\begin{align*}
\lim_{x \to -\infty} \limsup_{n \to \infty} \mu^{p^{(n)}} \left(\sup_{y\leq [nx]}A_iS_{y} \ge 0 \right)
&\le \lim_{x \to -\infty} \limsup_{n \to \infty}  \frac{2n\kappa}{[nx](c_0-c_i)^2}\left(\frac{2}{\kappa+1}+\frac{c_0+c_i}{\sqrt {n\kappa}}\right) \\
={}&\lim_{x \to -\infty} \frac{2\kappa}{x(c_0-c_i)^2}\frac{2}{\kappa+1}=0.
\end{align*}
\end{proof}

\begin{proof}[Proof of Theorem \ref{BMinvariant}]
Since $\mu^{p^{(n)}}$ is invariant under $T_i$, so is $\nu_n=\mu^{p^{(n)}}_{\frac{\sqrt\kappa}{\sqrt n},n}$  by Lemma \ref{ab}. Then Theorem \ref{donsker}, Lemma \ref{generaltheory} and \ref{checkassumption} show $\nu_{D}$ is also invariant under $T_i$.
In other words, two-sided standard $\kappa-$dimensional Brownian motion with drift $D\in\mathcal{D}$, given by
\[D=c_0e_0+\cdots+c_\kappa e_\kappa,\ c_0>c_i,\ \forall i\in\mathcal{C},\ c_0+\cdots+c_\kappa=0\]
is invariant under $T_i$. 
\end{proof}

\section*{Acknowledgements}

The author appreciates M.Sasada for her guidance and constructive comments from beggining to end. He also thanks D.Croydon for useful advice.

\vspace{10pt}

\end{document}